\documentclass[hidelinks,onefignum,onetabnum]{siamart251216}

\usepackage{lipsum}
\usepackage{amsfonts}
\usepackage{graphicx}
\usepackage{epstopdf}
\usepackage{algorithmic}
\usepackage{subcaption}
\usepackage{booktabs}
\ifpdf
  \DeclareGraphicsExtensions{.eps,.pdf,.png,.jpg}
\else
  \DeclareGraphicsExtensions{.eps}
\fi

\newsiamremark{remark}{Remark}
\newsiamremark{hypothesis}{Hypothesis}
\crefname{hypothesis}{Hypothesis}{Hypotheses}
\newsiamthm{claim}{Claim}
\newsiamremark{fact}{Fact}
\crefname{fact}{Fact}{Facts}

\headers{Hyper WENO5 for Conservation Law}{Y. Chen, W. Guo, and X. Zhong}

\title{Hypernetwork-Conditioned WENO5 Conservative-Form CNNs for One-Dimensional Conservation Laws\thanks{Submitted to the editors May 2026.}}

\author{
Yongsheng Chen
\thanks{School of Mathematical Sciences, Zhejiang University, Hangzhou, 310027, China. (\email{22035024@zju.edu.cn})}
\and
Wei Guo
\thanks{Corresponding author. Department of Mathematics and Statistics, Texas Tech University, Lubbock, TX, 70409, USA. (\email{weimath.guo@ttu.edu})}
\and
Xinghui Zhong
\thanks{School of Mathematical Sciences, Zhejiang University, Hangzhou, 310027, China. (\email{zhongxh@zju.edu.cn})}
}

\usepackage{amsopn}

\ifpdf
\hypersetup{
  pdftitle={Hyper WENO5 for Conservation Law},
  pdfauthor={Y. Chen, W. Guo, and X. Zhong}
}
\fi

\begin{document}

\maketitle

\begin{abstract}
We study a conservative data-driven discretization for one-dimensional hyperbolic conservation laws based on the classical fifth-order WENO finite-volume scheme and a hypernetwork architecture. In the proposed Hyper--WENO5 Conservative-Form Convolutional Neural Network (Hyper--CFCNN), a lightweight target network predicts the nonlinear WENO weights on each stencil, while a hypernetwork generates the target-network parameters from problem metadata, including the mesh spacing, mesh layout, and coarse descriptors of the initial condition. The construction preserves the standard polynomial reconstruction and conservative flux-difference update of WENO, which enables adaptation across problem instances and spatial resolutions without retraining. We also consider an unknown-flux variant, Hyper--CFCNN--F, in which a compact FluxNet is used in place of the analytical flux inside the numerical flux function while retaining a conservative update form. To improve long-time prediction quality, training uses a multi-step recurrent loss that penalizes error accumulation over successive time advances. Numerical experiments on one-dimensional test problems, including single- and multi-shock Burgers equations, the shallow-water system, and the Shu--Osher Euler example, show that Hyper--CFCNN attains accuracy comparable to classical WENO5, achieves near machine-precision conservation in the known-flux setting on fine meshes, and generalizes to unseen spatial resolutions and initial conditions without retraining. The flux-learning variant remains stable on meshes outside the training set and exhibits bounded conservation drift. These results show that hypernetwork-conditioned conservative WENO discretizations provide an effective framework for adaptive high-order learning of nonlinear conservation laws with either known or unknown fluxes.
\end{abstract}



\section{Introduction}

Hyperbolic conservation laws are a basic class of partial differential equations
(PDEs) arising in fluid dynamics, gas dynamics, traffic flow, plasma physics,
and related applications
\cite{leveque2002finite,godlewski2013numerical}. Even when the flux function is
smooth, smooth initial data may develop shocks, contact discontinuities, and
rarefaction waves in finite time \cite{lax1957hyperbolic}. This behavior makes
accurate numerical approximation difficult: high-order schemes can generate
nonphysical oscillations near discontinuities, whereas overly diffusive schemes
smear sharp solution features \cite{gottlieb1998total,lax1954weak}.

The design of accurate and stable discretizations for conservation laws has
therefore remained a central topic in scientific computing. High-resolution
methods, including total variation diminishing (TVD) schemes
\cite{harten1997high} and essentially non-oscillatory (ENO) reconstructions
\cite{shu1988efficient}, improve the balance between accuracy and robustness.
Among these methods, weighted essentially non-oscillatory (WENO) schemes are
widely used. Originally introduced by Liu et al.~\cite{liu1994weighted} and
refined by Jiang and Shu \cite{jiang1996efficient}, WENO methods construct
high-order interface values through nonlinear convex combinations of candidate
polynomial reconstructions on multiple stencils. In particular, the fifth-order
variant WENO5 is widely used because of its favorable accuracy--cost trade-off
\cite{shu2006essentially}. WENO schemes have been applied in a broad range of
problems, including compressible flows \cite{toth1996comparison},
magnetohydrodynamics \cite{balsara2000monotonicity}, relativistic flows
\cite{mignone2005hllc}, and shallow-water systems \cite{xing2005high}.

Classical WENO schemes, however, rely on prescribed analytical formulae for
smoothness indicators and nonlinear reconstruction weights. These constructions
are effective and well understood, but they do not directly incorporate
problem-dependent contextual information such as mesh characteristics or coarse
descriptors of the solution. In addition, the evaluation of smoothness indicators
and nonlinear weights contributes nontrivially to the cost of high-order
reconstruction, especially in more complex settings
\cite{suresh1997accurate,henrick2005mapped,borges2008improved}. Improved
variants such as mapped WENO and WENO-Z modify the analytical weighting
mechanism, but they remain hand-designed extensions of the same general
framework \cite{henrick2005mapped,borges2008improved,kossaczka2021enhanced}.

Recent machine-learning methods have introduced new tools for PDE modeling
and computation. Physics-informed neural networks (PINNs)
\cite{raissi2019physics} incorporate PDE residuals into training objectives and
have been used in both forward and inverse problems
\cite{karniadakis2021physics,yu2022gradient}. Neural operator models, including
Fourier neural operators (FNO) \cite{li2021fourier}, DeepONet
\cite{lu2021learning}, and graph neural operators (GNO)
\cite{sharma2024graph}, aim to approximate solution operators directly from
data. These approaches are useful in many settings, but they do not
automatically inherit the conservative structure and numerical stability
mechanisms built into classical finite-volume discretizations. For this reason, a
complementary direction is to integrate learning directly into
structure-preserving numerical schemes, while retaining the conservative update
form of the underlying discretization.

{\sloppy Recent work in this direction includes learned flux functions
\cite{morand2024deep}, neural-network-enhanced reconstruction procedures
\cite{li2023six}, and learning-based adaptive mesh refinement
\cite{foucart2023deep}. The DeepWENO framework
\cite{kossaczka2021enhanced}, for example, replaces the classical
smoothness-indicator pipeline with neural predictors while retaining the WENO
reconstruction template. Related learned finite-volume approaches have also been
investigated for improving reconstruction quality or reducing computational cost
\cite{stevens2020enhancement,ranade2021discretizationnet,chen2023learned}.
A central difficulty for such learned discretizations is transfer across problem
instances and spatial resolutions without retraining. \par}

Hypernetworks provide a natural mechanism for this purpose by generating the
parameters of a target network from contextual inputs
\cite{ha2017hypernetworks}. They have been studied in few-shot learning
\cite{munkhdalai2017meta}, multi-task learning
\cite{mahabadi2021parameter}, and continual learning
\cite{Oswald2020Continual}, and more recently have been explored in scientific
computing for parametric PDE solvers \cite{cho2023hypernetwork,de2021hyperpinn},
neural ODE models \cite{choromanski2020ode}, and adaptive discretizations in
fluid dynamics \cite{majumdar2023pihlora,chen2022meta}. These developments
suggest that conditioning numerical components on problem metadata may provide
a practical way to accommodate variability across meshes, physical parameters,
and initial conditions.

Motivated by these developments, we study a conservative, data-driven extension
of the WENO5 finite-volume scheme. We first consider a WENO5
conservative-form convolutional neural network (WENO5--CFCNN), in which a
lightweight neural predictor replaces the classical map from smoothness
indicators to nonlinear reconstruction weights, while the candidate WENO
polynomials and finite-volume flux-difference update are left unchanged. We
then introduce Hyper--CFCNN, in which a hypernetwork generates the parameters
of the weight-prediction network from problem metadata, including mesh
characteristics and descriptors of the initial state. This design is intended to
support transfer across problem instances and spatial resolutions without
retraining. We also study an unknown-flux variant, Hyper--CFCNN--F, in which a
compact neural model replaces analytical flux evaluation inside the numerical
flux computation, targeting settings where the governing flux is unavailable,
expensive to evaluate, or accessible only through data
\cite{ranade2021discretizationnet,morand2024deep}. In all variants, the update
is written in telescoping conservative flux-difference form; in the learned-flux
setting, discrete conservation is maintained when neighboring cells share the
same interface flux.

The paper examines the conservation properties, accuracy, and transfer behavior
of these learned discretizations on one-dimensional benchmark problems. The
numerical results show that hypernetwork-conditioned parameter generation can
be combined with conservative WENO-type updates while preserving the
underlying finite-volume structure.

The remainder of the paper is organized as follows. We begin in Section~\ref{sec:algorithm} with the necessary background, surveying hypernetwork architectures and finite-volume discretizations for conservation laws. Building on this foundation, Section~\ref{sec:hyper-weno5} introduces the proposed methods, covering the network architectures, training procedures, and boundary treatments in detail. The subsequent experimental evaluation in Section~\ref{sec:experiments} examines accuracy, conservation diagnostics, and generalization across a range of benchmark problems. Finally, Section~\ref{sec:con} summarizes our findings and outlines directions for future work.

\section{Background}
\label{sec:algorithm}

\subsection{Hypernetworks}

Hypernetworks are meta-learning architectures in which one network (the
\emph{hypernetwork}) generates the parameters of another network (the
\emph{target network}) conditioned on contextual information
\cite{ha2017hypernetworks}. In the present work, this mechanism is used to
generate the parameters of a lightweight target module within the learned
conservative discretization introduced in Section~\ref{sec:hyper-weno5}.

\begin{figure}[!htbp]
\centering
\includegraphics[width=0.8\textwidth]{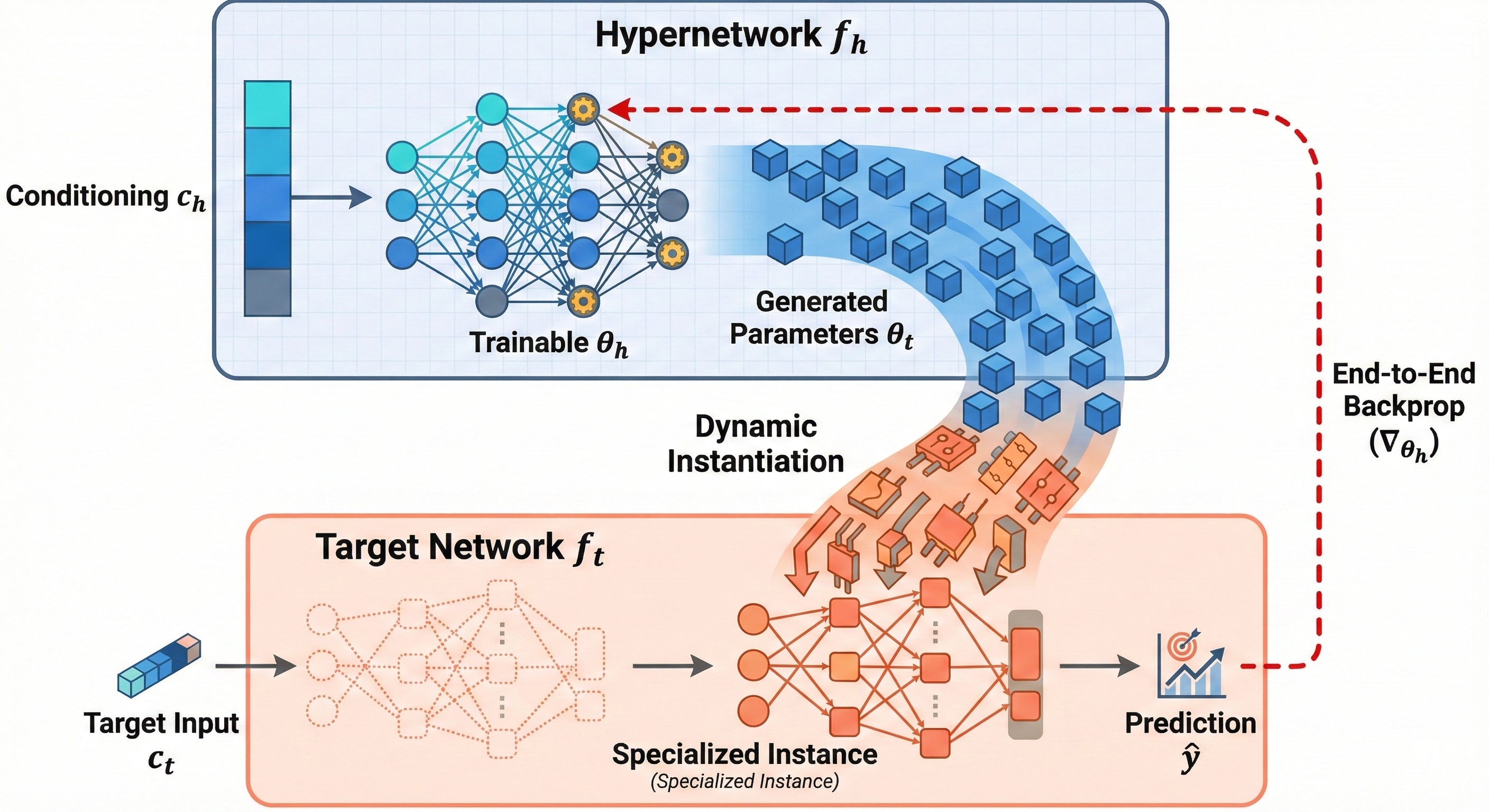}
\caption{Schematic of a hypernetwork. The hypernetwork \(f_h\) maps the
conditioning input \(c_h\) to the parameter vector
\(\boldsymbol{\theta}_t\), which instantiates the target network \(f_t\). The
target network then maps its input \(c_t\) to the output \(\hat y\). This
diagram summarizes the two-stage mapping in \eqref{eq:hyper} and the role of
the conditioning information in parameter generation.}
\label{fig:hyper}
\end{figure}

Formally, let \(f_h:\mathbb{R}^{d_h}\to\mathbb{R}^{P}\) denote the hypernetwork
and \(f_t:\mathbb{R}^{d_t}\to\mathbb{R}^{d_o}\) the target network, where
\(P\) is the number of parameters in \(f_t\). Given a conditioning vector
\(c_h\in\mathbb{R}^{d_h}\), the hypernetwork produces a parameter vector
\(\boldsymbol{\theta}_t\in\mathbb{R}^{P}\) for the target network, which then
maps an input \(c_t\in\mathbb{R}^{d_t}\) to an output
\(\hat y\in\mathbb{R}^{d_o}\):
\begin{equation}\label{eq:hyper}
\boldsymbol{\theta}_t = f_h(c_h;\boldsymbol{\theta}_h),
\qquad
\hat y = f_t(c_t;\boldsymbol{\theta}_t),
\end{equation}
where \(\boldsymbol{\theta}_h\) denotes the trainable parameters of the
hypernetwork. The corresponding composite mapping is
\[
(c_h,c_t)\mapsto f_t\bigl(c_t; f_h(c_h;\boldsymbol{\theta}_h)\bigr),
\]
which can be trained end-to-end by backpropagation with respect to
\(\boldsymbol{\theta}_h\), provided the constituent mappings are differentiable.
Figure~\ref{fig:hyper} illustrates this two-stage construction, in which the
conditioning input \(c_h\) determines the parameter vector
\(\boldsymbol{\theta}_t\) and thereby instantiates the target network used on
the input \(c_t\).

In this construction, the conditioning variable \(c_h\) modulates the parameters
of the target network and thereby enables instance-dependent adaptation. 

\subsection{Conservation laws and the WENO5 scheme}
\label{sec:background-weno}

{\sloppy For simplicity, we present the scalar one-dimensional conservation law
\begin{equation}\label{eq:conservation}
u_t + \bigl(f(u)\bigr)_x = 0, \qquad x\in\Omega,\ t>0,
\end{equation}
supplemented with appropriate initial and boundary conditions. Let
\(I_i=[x_{i-\frac12},x_{i+\frac12}]\), \(i=1,\dots,N\), denote a uniform
partition of the computational domain, with cell centers \(x_i\) and mesh size
\(\Delta x\). The corresponding cell average is \par}
\begin{equation}\label{eq:cell_avg}
\bar u_i(t)=\frac{1}{\Delta x}\int_{I_i}u(x,t)\,dx.
\end{equation}

Integrating \eqref{eq:conservation} over \(I_i\) gives the semi-discrete
finite-volume scheme
\begin{equation}\label{eq:semi_discrete}
\frac{d\bar u_i}{dt}
=
-\frac{1}{\Delta x}\left(\hat f_{i+\frac12}-\hat f_{i-\frac12}\right),
\end{equation}
where \(\hat f_{i+\frac12}\) is a consistent numerical flux at the interface
\(x_{i+\frac12}\). Summing \eqref{eq:semi_discrete} over all cells produces
telescoping flux differences. Under periodic or homogeneous no-flux (zero-gradient) boundary
conditions, the discrete total mass \(\sum_i \bar u_i\) is therefore conserved,
provided that each shared interface is assigned a single numerical flux.

To obtain high-order accuracy, one reconstructs one-sided interface states from
cell averages. In WENO5, the left-biased state at \(x_{i+\frac12}\) is written as
a nonlinear convex combination of three third-order candidate polynomials on
the stencils \(S_k=\{I_{i+k-2},I_{i+k-1},I_{i+k}\}\), \(k=0,1,2\):
\begin{equation}\label{eq:weno_reconstruct}
u^-_{i+\frac12}
=
\sum_{k=0}^{2}\omega_k^-\, q_k^-\!\left(x_{i+\frac12}\right).
\end{equation}
For cells adjacent to a physical boundary, these stencils extend outside the
computational domain. In that case, ghost-cell values are supplied in a
boundary-consistent manner so that the reconstruction remains well defined and
retains its designed order. The boundary treatment used in the proposed method
is described in Section~\ref{sec:bc-fluxnet}.

On a uniform mesh, the candidate polynomials evaluated at \(x_{i+\frac12}\) are
\begin{align}
q_0^-(x_{i+\frac12}) &= \tfrac{1}{6}\bigl(2\bar{u}_{i-2}-7\bar{u}_{i-1}+11\bar{u}_i\bigr), \label{eq:poly0}\\
q_1^-(x_{i+\frac12}) &= \tfrac{1}{6}\bigl(-\bar{u}_{i-1}+5\bar{u}_i+2\bar{u}_{i+1}\bigr), \label{eq:poly1}\\
q_2^-(x_{i+\frac12}) &= \tfrac{1}{6}\bigl(2\bar{u}_i+5\bar{u}_{i+1}-\bar{u}_{i+2}\bigr). \label{eq:poly2}
\end{align}
The nonlinear weights are defined by
\begin{equation}\label{eq:smoothness_indicators}
\omega_k^-=\frac{\tilde\omega_k^-}{\sum_{j=0}^{2}\tilde\omega_j^-},
\qquad
\tilde\omega_k^-=\frac{d_k}{(\epsilon+\beta_k^-)^r},
\end{equation}
where \(d_k\) are the optimal linear weights, \(\beta_k^-\) are the smoothness
indicators, \(\epsilon>0\) is a small regularization parameter, and \(r\)
controls the nonlinearity of the weighting. In standard WENO5, one typically
takes \(\epsilon=10^{-6}\) and \(r=2\). In smooth regions, this construction
recovers fifth-order accuracy, while near nonsmooth features the nonlinear
weights bias the reconstruction toward smoother substencils
\cite{shu2006essentially}.

The right-biased state \(u^+_{i+\frac12}\) is constructed analogously using the
mirrored stencils
\[
S_k^{+}=\{I_{i+1-k},\,I_{i+2-k},\,I_{i+3-k}\},\qquad k=0,1,2,
\]
with candidate polynomials
\begin{align}
q_0^+(x_{i+\frac12}) &= \tfrac{1}{6}\bigl(11\bar{u}_{i+1}-7\bar{u}_{i+2}+2\bar{u}_{i+3}\bigr), \\
 q_1^+(x_{i+\frac12}) &= \tfrac{1}{6}\bigl(2\bar{u}_{i}+5\bar{u}_{i+1}-\bar{u}_{i+2}\bigr), \\
 q_2^+(x_{i+\frac12}) &= \tfrac{1}{6}\bigl(-\bar{u}_{i-1}+5\bar{u}_{i}+2\bar{u}_{i+1}\bigr).
\end{align}
The corresponding weights \(\omega_k^+\) are defined analogously using the
optimal linear weights and smoothness indicators associated with the mirrored
stencils. The right-biased reconstruction is then
\[
u^+_{i+\frac12}
=
\sum_{k=0}^{2}\omega_k^+\, q_k^+\!\left(x_{i+\frac12}\right).
\]

Given the reconstructed states \((u^-_{i+\frac12},u^+_{i+\frac12})\), we employ
the local Lax--Friedrichs (Rusanov) flux \cite{lax1954weak}
\begin{equation}\label{eq:lax_friedrichs}
\hat f_{i+\frac12}
=
\frac12\left[
f\!\left(u^-_{i+\frac12}\right)
+
f\!\left(u^+_{i+\frac12}\right)
-
\alpha_{i+\frac12}\left(u^+_{i+\frac12}-u^-_{i+\frac12}\right)
\right],
\end{equation}
where \(\alpha_{i+\frac12}\) is an upper bound on the characteristic speed, for
example,
\[
\alpha_{i+\frac12}
=
\max\Bigl(
|f'(u^-_{i+\frac12})|,
|f'(u^+_{i+\frac12})|
\Bigr).
\]
For systems, one typically replaces this quantity by the spectral radius of the
flux Jacobian.

The resulting semi-discrete system \eqref{eq:semi_discrete} is advanced in time
by the third-order strong-stability-preserving Runge--Kutta method (SSP--RK3)
\cite{shu1988efficient}:
\begin{align}
\bar{\mathbf{u}}^{(1)}
&=
\bar{\mathbf{u}}^n+\Delta t\,\mathcal{F}(\bar{\mathbf{u}}^n), \nonumber\\
\bar{\mathbf{u}}^{(2)}
&=
\tfrac{3}{4}\bar{\mathbf{u}}^n
+
\tfrac{1}{4}\Bigl(\bar{\mathbf{u}}^{(1)}+\Delta t\,\mathcal{F}(\bar{\mathbf{u}}^{(1)})\Bigr), \nonumber\\
\bar{\mathbf{u}}^{n+1}
&=
\tfrac{1}{3}\bar{\mathbf{u}}^n
+
\tfrac{2}{3}\Bigl(\bar{\mathbf{u}}^{(2)}+\Delta t\,\mathcal{F}(\bar{\mathbf{u}}^{(2)})\Bigr),
\label{eq:ssprk3}
\end{align}
where \(\bar{\mathbf{u}}=(\bar{u}_1,\ldots,\bar{u}_N)^\top\) and
\(\mathcal{F}\) denotes the spatial discretization operator in
\eqref{eq:semi_discrete}.

\section{Hyper--WENO5 Conservative--Form Convolutional Neural Network}
\label{sec:hyper-weno5}

\begin{figure}[!htbp]
\centering
\includegraphics[width=\textwidth]{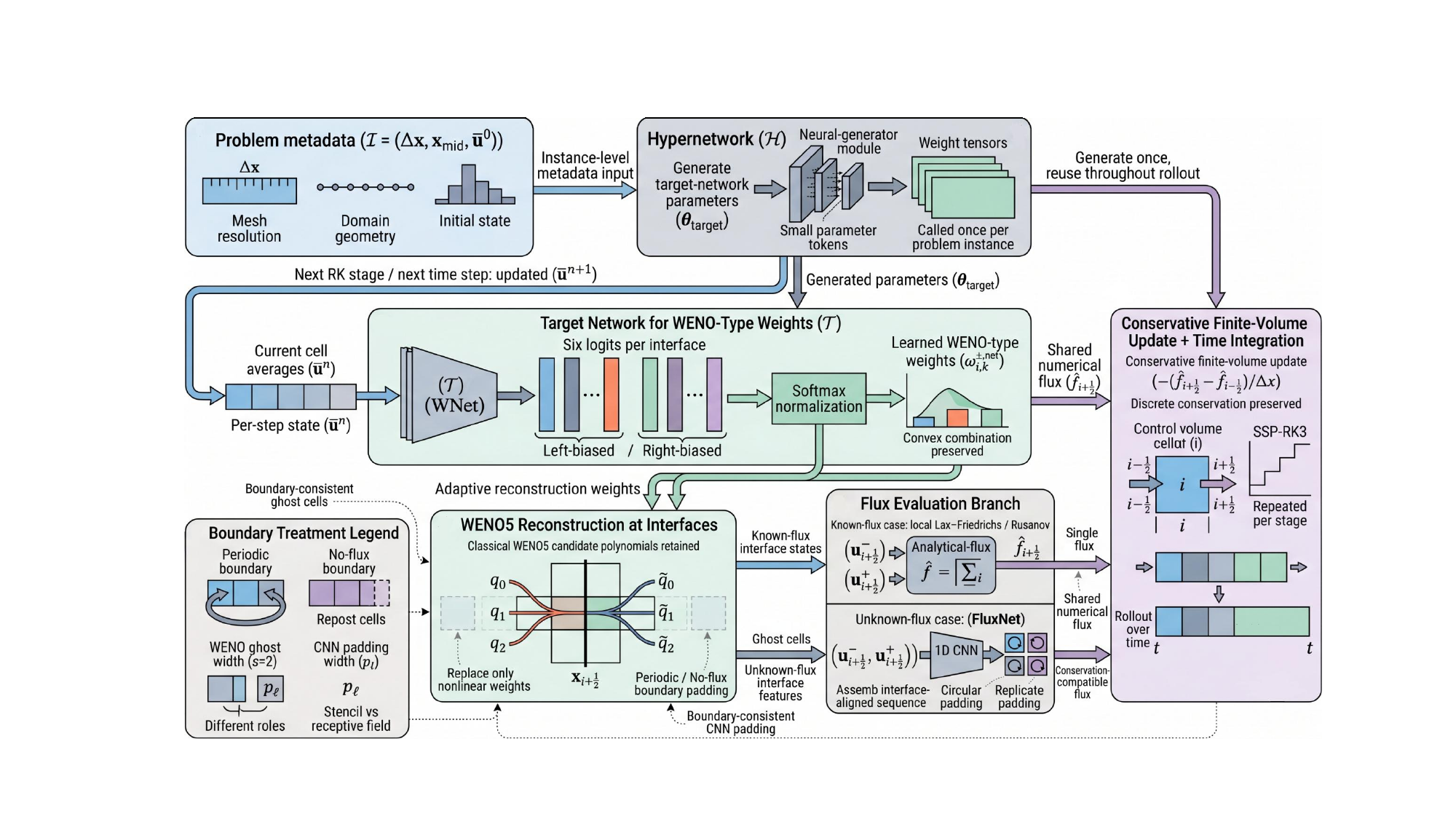}
\caption{Overall framework of Hyper--CFCNN and its unknown-flux variant
Hyper--CFCNN--F. The hypernetwork generates the parameters of a lightweight
target network from problem metadata, and the resulting reconstruction weights
are used in a conservative finite-volume update with either analytical or
learned interface fluxes.}
\label{fig:hyper_cfn}
\end{figure}

This section introduces the proposed learned conservative discretization
framework built on the finite-volume WENO5 formulation reviewed in
Section~\ref{sec:background-weno}. Figure~\ref{fig:hyper_cfn} gives an overview of the full framework, including the hypernetwork-conditioned known-flux model Hyper--CFCNN and its unknown-flux variant Hyper--CFCNN--F. The framework consists of three main
ingredients: a target network that predicts WENO-type reconstruction weights, a
hypernetwork that generates the parameters of this target network from problem
metadata, and a conservative finite-volume update that uses either analytical
or learned interface fluxes. We begin with the base reconstruction model
WENO5--CFCNN, then introduce the hypernetwork-conditioned parameterization, and
finally describe the unknown-flux extension, the forward pass with boundary
treatment, and the training procedure.

\subsection{WENO5 Conservative--Form Convolutional Neural Network}
\label{sec:weno5-cfcnn}

We first define a learned reconstruction model that serves as the base
component of the full framework. For clarity, we present the scalar case;
extensions to systems are carried out componentwise or in characteristic
variables in the standard manner.

Starting from the semi-discrete finite-volume formulation
\eqref{eq:semi_discrete}, we define a shallow one-dimensional convolutional
neural network \(\textsf{WNet}\) that acts on the cell-average vector
\(\bar{\boldsymbol{u}}=(\bar u_1,\dots,\bar u_N)^\top\). For each interface
\(x_{i+\frac12}\), the network outputs six logits corresponding to the left-
and right-biased WENO5 reconstruction weights. To avoid conflict with the
classical WENO notation \(\omega_k^\pm\), we denote the network output by
\begin{equation}\label{eq:wnet_logits}
\boldsymbol{\eta}_i
=
\big[\textsf{WNet}(\bar{\boldsymbol{u}};\boldsymbol{\theta}_{\mathrm{w}})\big]_i
\in \mathbb{R}^{6},
\qquad
\boldsymbol{\eta}_i
=
\big(
\eta_{i,0}^{-},\eta_{i,1}^{-},\eta_{i,2}^{-},
\eta_{i,0}^{+},\eta_{i,1}^{+},\eta_{i,2}^{+}
\big)^\top,
\end{equation}
where \(\boldsymbol{\theta}_{\mathrm{w}}\) denotes the trainable parameters of
\(\textsf{WNet}\).

The learned reconstruction weights are obtained by applying separate softmax
normalizations to the left- and right-biased logits:
\begin{equation}\label{eq:softmax_net}
\omega_{i,k}^{-,\mathrm{net}}
=
\frac{\exp(\eta_{i,k}^{-})}{\sum_{m=0}^{2}\exp(\eta_{i,m}^{-})},
\qquad
\omega_{i,k}^{+,\mathrm{net}}
=
\frac{\exp(\eta_{i,k}^{+})}{\sum_{m=0}^{2}\exp(\eta_{i,m}^{+})},
\qquad k=0,1,2.
\end{equation}
By construction,
\(\sum_{k=0}^{2}\omega_{i,k}^{\pm,\mathrm{net}}=1\) and
\(\omega_{i,k}^{\pm,\mathrm{net}}\ge 0\), so the convex-combination structure
of the WENO reconstruction is preserved.

Only the classical nonlinear weights are replaced. The candidate WENO5
polynomials \(q_k^-\) and \( q_k^+\) introduced in
Section~\ref{sec:background-weno} are retained unchanged. The reconstructed
left- and right-biased interface states are therefore
\begin{equation}\label{eq:left_rec_net}
u^-_{i+\frac12}
=
\sum_{k=0}^{2}\omega_{i,k}^{-,\mathrm{net}}\,q_k^-\!\bigl(x_{i+\frac12}\bigr),
\end{equation}
and
\begin{equation}\label{eq:right_rec_net}
u^+_{i+\frac12}
=
\sum_{k=0}^{2}\omega_{i,k}^{+,\mathrm{net}}\,q_k^+\!\bigl(x_{i+\frac12}\bigr).
\end{equation}

The spatial update remains in conservative finite-volume form:
\begin{equation}\label{eq:fv_net}
\frac{d\bar{u}_i}{dt}
=
-\frac{1}{\Delta x}\bigl(\hat f_{i+\frac12}-\hat f_{i-\frac12}\bigr).
\end{equation}
When the analytical flux \(f(\cdot)\) is available, we use the local
Lax--Friedrichs (Rusanov) flux
\begin{equation}\label{eq:lff_net}
\hat f_{i+\frac12}
=
\frac{1}{2}\left[
f\!\left(u^-_{i+\frac12}\right)
+
f\!\left(u^+_{i+\frac12}\right)
-
\alpha_{i+\frac12}\left(u^+_{i+\frac12}-u^-_{i+\frac12}\right)
\right],
\end{equation}
where, for scalar problems,
\[
\alpha_{i+\frac12}
=
\max\Bigl\{
|f'(u^-_{i+\frac12})|,
|f'(u^+_{i+\frac12})|
\Bigr\}.
\]
Time integration is performed by SSP--RK3 as in \eqref{eq:ssprk3}. Since the
update is kept in flux-difference form, the semi-discrete scheme remains
discretely conservative whenever a single shared numerical flux is assigned to
each interface.

{\sloppy This construction defines the basic learned reconstruction mechanism used throughout the paper. The next step is to replace the directly learned target-network parameters \(\boldsymbol{\theta}_{\mathrm{w}}\) by an instance-dependent parameter set generated from problem metadata. \par }

\subsection{Hypernetwork-conditioned weight prediction}
\label{sec:hyper-cfcnn}

{\sloppy We now extend WENO5--CFCNN by introducing a hypernetwork-conditioned
parameterization. The target network \(\mathcal T\) has the same architecture
as \(\textsf{WNet}\) in Section~\ref{sec:weno5-cfcnn}, but its parameters are no longer optimized directly. Instead, they are generated once
per problem instance by a hypernetwork \(\mathcal H\). This separates
instance-level parameter generation from per-step numerical evaluation.\par }

We encode the problem metadata as
\begin{equation}\label{eq:hyper_input}
\mathcal I
=
\bigl(\Delta x,\mathbf{x}_{\mathrm{mid}},\bar{\mathbf{u}}^{\,0}\bigr),
\end{equation}
where \(\mathbf{x}_{\mathrm{mid}}=(x_1,\ldots,x_N)^\top\) is the vector of cell
centers and \(\bar{\mathbf{u}}^{\,0}\) is the initial cell-average state on the
computational mesh. Thus, \(\mathcal I\) encodes the computational mesh together with the rollout-initial state.

The hypernetwork maps \(\mathcal I\) to the target-network parameters:
\begin{equation}\label{eq:hyper_out}
\boldsymbol{\theta}_{\mathrm{target}}
=
\mathcal H(\mathcal I;\boldsymbol{\theta}_h),
\end{equation}
where \(\boldsymbol{\theta}_h\) denotes the trainable parameters of the
hypernetwork. During time stepping, the target network then produces
per-interface logits
\begin{equation}\label{eq:target_call}
\boldsymbol{\eta}^{\,n}
=
\mathcal T(\bar{\mathbf{u}}^{\,n};\boldsymbol{\theta}_{\mathrm{target}})
\in\mathbb{R}^{N\times 6},
\end{equation}
which are converted to bias-wise normalized reconstruction weights through
\eqref{eq:softmax_net}. These weights are then used in
\eqref{eq:left_rec_net}--\eqref{eq:right_rec_net} and in the conservative
update \eqref{eq:fv_net}.

{\sloppy Thus, the hypernetwork is evaluated once per problem instance to generate
\(\boldsymbol{\theta}_{\mathrm{target}}\), whereas the lightweight target
network is evaluated repeatedly throughout time integration. The numerical
update remains local at each step, but the target-network parameters depend on
instance-level metadata. This yields the hypernetwork-conditioned formulation
Hyper--CFCNN. \par}

\subsection{Unknown-flux extension: Hyper--CFCNN--F}
\label{sec:fluxnet-subsub}

{\sloppy We now consider the case in which the analytical flux \(f(\cdot)\) is unavailable, expensive to evaluate, or accessible only through data. In this setting, we
define an unknown-flux variant of Hyper--CFCNN, denoted by Hyper--CFCNN--F, by
replacing the analytical flux evaluation in \eqref{eq:lff_net} with a learned
interface-flux model:\par }
\begin{equation}\label{eq:fluxnet}
\hat f_{i+\frac12}
=
\textsf{FluxNet}\!\left(
u^-_{i+\frac12},u^+_{i+\frac12};
\boldsymbol{\theta}_{\mathrm{flux}}
\right),
\end{equation}
where \(\boldsymbol{\theta}_{\mathrm{flux}}\) are the trainable parameters of
\(\textsf{FluxNet}\). The finite-volume update remains the conservative
flux-difference form \eqref{eq:fv_net}. Accordingly, discrete conservation is
preserved provided that each interface \(x_{i+\frac12}\) is assigned a single
shared numerical flux used consistently by the two neighboring cells.

{ \sloppy The known-flux and unknown-flux variants therefore share the same reconstruction
and finite-volume update structure. Their difference lies only in the
flux-evaluation branch. \par}

\subsection{Forward pass and boundary treatment}
\label{sec:forward}

We now summarize the rollout procedure shared by Hyper--CFCNN and
Hyper--CFCNN--F. At the beginning of each rollout, the hypernetwork generates
the target-network parameters once for the given problem instance:
\[
\boldsymbol{\theta}_{\mathrm{target}}
=
\mathcal H(\mathcal I;\boldsymbol{\theta}_h).
\]
These parameters are then reused throughout all SSP--RK3 stage evaluations in
that rollout.

At each Runge--Kutta stage, the forward pass consists of four steps. First, the
target network is evaluated on the current cell averages,
\[
\boldsymbol{\eta}^{\,n}
=
\mathcal T(\bar{\mathbf{u}}^{\,n};\boldsymbol{\theta}_{\mathrm{target}}),
\]
and the normalized learned weights
\(\omega_{i,k}^{\pm,\mathrm{net}}\) are obtained through
\eqref{eq:softmax_net}. Second, boundary-consistent ghost cells are introduced,
and the interface states \(u^-_{i+\frac12}\) and \(u^+_{i+\frac12}\) are
reconstructed using \eqref{eq:left_rec_net}--\eqref{eq:right_rec_net}. Third,
the interface fluxes are evaluated, using \eqref{eq:lff_net} in the
known-flux case and \eqref{eq:fluxnet} in the unknown-flux case. Fourth, the cell averages are updated through the conservative finite-volume form \eqref{eq:fv_net}. The two variants differ only in the third step. The second and third steps require boundary treatment for the WENO reconstruction and, in the unknown-flux setting, for the CNN-based flux evaluation.


\subsubsection{Boundary conditions and \textsf{FluxNet} design}
\label{sec:bc-fluxnet}

Boundary treatment enters the method in two distinct ways. First, the WENO5
reconstruction requires ghost cells because the reconstruction stencil extends
beyond the physical domain near the boundary. Second, in the unknown-flux
setting, the convolutional layers in \textsf{FluxNet} require padding so that
their receptive fields remain compatible with the prescribed boundary
condition. These two mechanisms serve different purposes and are therefore
specified separately below.

\paragraph{Notation}
Let \(s\) denote the half-width of the WENO reconstruction stencil. For WENO5,
\(s=2\). Ghost cells for the finite-volume variables are padded with width \(s\)
on each side. Separately, let \(p_\ell\) denote the padding width in the
\(\ell\)-th convolutional layer of \(\textsf{FluxNet}\). For a one-dimensional
convolution with odd kernel size \(k_\ell\), unit stride, and no dilation, we
take
\[
p_\ell=\frac{k_\ell-1}{2}.
\]
Thus, \(s\) is determined by the discretization stencil, whereas \(p_\ell\) is
determined by the CNN architecture.

\paragraph{WENO5 reconstruction}
For WENO5, we pad \(s=2\) ghost cells on each side of the finite-volume state.
For systems, this is done componentwise. The ghost-cell rule depends on the
boundary condition.

\emph{Periodic boundaries.}
For \(j=1,\ldots,s\),
\begin{equation}\label{eq:periodic-pad}
\bar{u}_{1-j}=\bar{u}_{N+1-j},
\qquad
\bar{u}_{N+j}=\bar{u}_{j}.
\end{equation}

\emph{No-flux boundaries.}
For \(j=1,\ldots,s\),
\begin{equation}\label{eq:noflux-pad}
\bar{u}_{1-j}=\bar{u}_{1},
\qquad
\bar{u}_{N+j}=\bar{u}_{N}.
\end{equation}

\paragraph{\textsf{FluxNet}}
{\sloppy In the unknown-flux setting, \(\textsf{FluxNet}\) is implemented as a
one-dimensional CNN that maps interface-aligned features to interface-aligned
numerical fluxes. For scalar conservation laws, we define \par}
\[
X
=
\big\{
(u^-_{i+\frac12},u^+_{i+\frac12})
\big\}_{i=1}^{N}
\in\mathbb{R}^{N\times 2}.
\]
The predicted interface-flux sequence is
\begin{equation}\label{eq:fluxnet_cnn}
\hat{\boldsymbol{f}}
=
\textsf{FluxNet}(X;\boldsymbol{\theta}_{\mathrm{flux}})
\in\mathbb{R}^{N},
\qquad
\hat f_{i+\frac12}
=
[\hat{\boldsymbol{f}}]_i.
\end{equation}

Although \(X\) is assembled from WENO reconstructions that already use ghost
cells, each convolutional layer in \(\textsf{FluxNet}\) still requires padding so
that the output length remains \(N\) and the receptive field remains compatible
with the boundary condition.

\begin{itemize}
\item \emph{Periodic boundaries.}
Use circular padding in every convolutional layer. For layer \(\ell\) with
kernel size \(k_\ell\), set \(p_\ell=(k_\ell-1)/2\). For
\(X\in\mathbb{R}^{N\times 2}\), the padded input is
\[
\widetilde{X}
=
\big[
X_{N-p_\ell+1:N}\,\big|\,X_{1:N}\,\big|\,X_{1:p_\ell}
\big].
\]

\item \emph{No-flux boundaries.}
Use replicate padding of width \(p_\ell\):
\[
\widetilde{X}
=
\big[
\underbrace{X_1,\ldots,X_1}_{p_\ell}
\,\big|\,
X_{1:N}
\,\big|\,
\underbrace{X_N,\ldots,X_N}_{p_\ell}
\big].
\]
This choice aligns the CNN boundary treatment with the no-flux boundary.
\end{itemize}

\paragraph{Implementation remarks}
Discrete conservation follows from the flux-difference form \eqref{eq:fv_net}.
The ghost-cell rules \eqref{eq:periodic-pad}--\eqref{eq:noflux-pad} affect only
how boundary information enters the WENO reconstruction, whereas the CNN-layer
padding described above affects only how boundary information enters the
\textsf{FluxNet} evaluation. Neither modification changes the conservative
flux-difference update itself. The stencil width \(s\) is determined by the
discretization, whereas the padding widths \(p_\ell\) are determined by the
network architecture. For systems, the ghost-cell rules
\eqref{eq:periodic-pad}--\eqref{eq:noflux-pad} are applied componentwise to the
finite-volume variables, while the CNN padding is applied along the spatial
dimension to each channel of the \textsf{FluxNet} input. If characteristic reconstruction is used, the same ghost-cell and CNN-padding rules are applied to the characteristic variables.
\subsection{Training methodology}
\label{sec:training}

We train Hyper--CFCNN by supervised trajectory matching on randomly windowed
segments extracted from reference solution trajectories, with the loss defined
on multi-step unrolled predictions. The trainable variables are the
hypernetwork parameters \(\boldsymbol{\theta}_h\) and, in the unknown-flux
setting, the FluxNet parameters \(\boldsymbol{\theta}_{\mathrm{flux}}\). The
target-network parameters are never optimized directly; they are generated by
the hypernetwork for each problem instance and mesh.

\paragraph{Multi-resolution reference data and windowing}
Let \(l\in\{1,\dots,L_{\mathrm{mesh}}\}\) index the mesh level, with
\(N_l\) cells, spacing \(\Delta x_l\), and cell centers
\(\mathbf{x}_{\mathrm{mid}}^{(l)}\in\mathbb{R}^{N_l}\). For each
training instance \(i=1,\dots,N_{\mathrm{traj}}\) and each mesh level \(l\), we
precompute a reference trajectory of cell averages using a high-fidelity
solver:
\[
\bigl\{
\bar{\mathbf{u}}^{(i,l)}(t_m^{(l)})
\bigr\}_{m=0}^{M_l},
\qquad
t_m^{(l)}=m\,\Delta t_l,
\]
where \(\Delta t_l\) satisfies the prescribed CFL condition on mesh level \(l\).
During training, we sample a start index
\(s\in\{0,\dots,M_l-L\}\) uniformly at random and form the window
\begin{equation}\label{eq:train_window}
\mathcal{W}_{i,l,s}
=
\Bigl\{
\bar{\mathbf{u}}^{(i,l)}(t_s^{(l)}),
\bar{\mathbf{u}}^{(i,l)}(t_{s+1}^{(l)}),
\ldots,
\bar{\mathbf{u}}^{(i,l)}(t_{s+L}^{(l)})
\Bigr\}.
\end{equation}
This windowing procedure exposes the model to local trajectory segments from
multiple mesh levels while controlling the computational cost of training.

\paragraph{One-step map and unrolling}
For the window \(\mathcal{W}_{i,l,s}\), we define the metadata
\begin{equation}\label{eq:meta_train}
\mathcal{I}_{i,l}
=
\bigl(
\Delta x_l,
\mathbf{x}_{\mathrm{mid}}^{(l)},
\bar{\mathbf{u}}^{(i,l)}(t_s^{(l)})
\bigr),
\end{equation}
which is held fixed over the corresponding rollout. Thus, the conditioning
metadata encode the mesh and the rollout-initial state, whereas the
instantaneous solution state evolves under the one-step update. The
hypernetwork then generates the corresponding target-network parameters
\[
\boldsymbol{\theta}_{\mathrm{target}}^{(i,l)}
=
\mathcal H(\mathcal I_{i,l};\boldsymbol{\theta}_h).
\]

We define the mesh-dependent one-step operator
\begin{equation}\label{eq:onestep}
\hat{\mathbf{u}}_{s+1}^{(i,l)}
=
\Phi_\Theta^{[l]}
\bigl(
\bar{\mathbf{u}}^{(i,l)}(t_s^{(l)});\mathcal I_{i,l}
\bigr),
\end{equation}
where
\[
\Theta
=
\begin{cases}
(\boldsymbol{\theta}_h), & \text{known-flux case},\\[2pt]
(\boldsymbol{\theta}_h,\boldsymbol{\theta}_{\mathrm{flux}}), & \text{unknown-flux case}.
\end{cases}
\]
Here \(\Phi_\Theta^{[l]}\) denotes one SSP--RK3 step on mesh level \(l\), using
the target-network parameters generated from \(\mathcal I_{i,l}\), the
reconstruction formulas \eqref{eq:left_rec_net}--\eqref{eq:right_rec_net}, and
either the analytical-flux evaluation \eqref{eq:lff_net} or the learned-flux
evaluation \eqref{eq:fluxnet_cnn}.

With initialization
\[
\hat{\mathbf{u}}_{s}^{(i,l)}
=
\bar{\mathbf{u}}^{(i,l)}(t_s^{(l)}),
\]
we unroll this one-step map for \(K\) steps:
\begin{equation}\label{eq:unroll_pred}
\hat{\mathbf{u}}_{s+\ell}^{(i,l)}
=
\bigl(\Phi_\Theta^{[l]}\bigr)^{(\ell)}
\bigl(
\bar{\mathbf{u}}^{(i,l)}(t_s^{(l)});\mathcal I_{i,l}
\bigr),
\qquad
\ell=1,\dots,K,
\end{equation}
where \((\cdot)^{(\ell)}\) denotes \(\ell\)-fold composition. In this work, we
take \(K=4\), with \(1\le K\le L\).

\paragraph{Multi-resolution recurrent loss}
To compare rollout errors across different meshes, we use a mesh-weighted
discrete \(L^2\) norm. For the prediction error at step \(\ell\), we write
\[
\left\|
\hat{\mathbf{u}}_{s+\ell}^{(i,l)}
-
\bar{\mathbf{u}}^{(i,l)}(t_{s+\ell}^{(l)})
\right\|_{2,l}^{2}
:=
\Delta x_l\sum_{j=1}^{N_l}
\left(
\hat{u}_{s+\ell,j}^{(i,l)}
-
\bar{u}_{j}^{(i,l)}(t_{s+\ell}^{(l)})
\right)^2,
\]
which provides a consistent approximation of the continuous \(L^2\) error
across mesh levels. The training objective is
\begin{equation}\label{eq:traj_loss}
\mathcal{L}_{\mathrm{traj}}(\Theta)
=
\frac{1}{N_{\mathrm{traj}}L_{\mathrm{mesh}}}
\sum_{i=1}^{N_{\mathrm{traj}}}
\sum_{l=1}^{L_{\mathrm{mesh}}}
\frac{1}{|\mathcal{S}_{i,l}|}
\sum_{s\in\mathcal{S}_{i,l}}
\frac{1}{K}
\sum_{\ell=1}^{K}
\left\|
\hat{\mathbf{u}}_{s+\ell}^{(i,l)}
-
\bar{\mathbf{u}}^{(i,l)}(t_{s+\ell}^{(l)})
\right\|_{2,l}^{2},
\end{equation}
where \(\mathcal{S}_{i,l}\subset\{0,\dots,M_l-L\}\) denotes the set of sampled
window starts for instance \(i\) on mesh level \(l\) within the current epoch.

\paragraph{Optimization and batching}
Mini-batches are formed by sampling tuples \((i,l,s)\) across trajectories and
mesh levels. Gradients are computed by backpropagation through the unrolled map
\((\Phi_\Theta^{[l]})^{(\ell)}\). At inference time, the same rollout
procedure is used: given metadata \(\mathcal I_{i,l}\), the hypernetwork emits
\(\boldsymbol{\theta}_{\mathrm{target}}^{(i,l)}\) once, after which the target
network is evaluated repeatedly at each Runge--Kutta stage on mesh level \(l\).
\section{Numerical results}\label{sec:experiments}

We evaluate H--CFCNN and its unknown-flux variant H--CFCNN--F on three classes
of one-dimensional conservation laws: the inviscid Burgers equation, the
shallow-water system, and the compressible Euler equations. The experiments
examine accuracy under mesh refinement, discrete conservation, transfer to
unseen meshes and initial conditions, and sensitivity to rollout horizon and
time step.

Unless stated otherwise, training follows the procedure in
Section~\ref{sec:training}, and inference uses the forward pass and boundary
treatment described in Section~\ref{sec:forward} and
Section~\ref{sec:bc-fluxnet}. For each benchmark, we specify the training data,
mesh hierarchy, rollout horizon, and reference solver before discussing the
numerical results.

\paragraph{Discrete conservation diagnostic}
Let \(\bar q_j(t)\) denote the cell average of a conserved quantity \(q\) on the
cell \(I_j\). Let \(\hat F_{1/2}^{(m)}\) and \(\hat F_{N+1/2}^{(m)}\) denote the
numerical boundary fluxes at time level \(t_m\), computed in the same manner as
the interior interface fluxes: by \eqref{eq:lff_net} in the known-flux setting
and by \eqref{eq:fluxnet_cnn} in the learned-flux setting. We measure discrete
conservation by
\begin{equation}\label{eq:disc_cons_remainder}
C(q)(t_\ell)
:=
\left|
\sum_{j=1}^N \big(\bar q_j(t_\ell)-\bar q_j(t_0)\big)\,\Delta x
-
\sum_{k=1}^{\ell}
\Big(\hat F_{1/2}^{(k-1)}-\hat F_{N+1/2}^{(k-1)}\Big)\,\Delta t
\right|.
\end{equation}
For periodic boundaries, the boundary-flux term vanishes identically. For no-flux boundaries, the boundary term is computed from the numerical fluxes actually used by the scheme, so any nonzero remainder measures the conservation error accumulated during the rollout.

\paragraph{Reference data and network components}
For each PDE family, reference trajectories are generated by a high-fidelity
finite-volume solver on a problem-dependent mesh hierarchy under a fixed CFL
condition. Training uses the multi-resolution windowing strategy and recurrent
rollout loss introduced in Section~\ref{sec:training}. The target network
\(\mathcal T\) outputs six WENO-type logits per cell, which are converted to
bias-wise normalized reconstruction weights by \eqref{eq:softmax_net}; its
parameter set is generated once per problem instance by the hypernetwork
\(\boldsymbol{\theta}_{\mathrm{target}}=\mathcal H(\mathcal I;
\boldsymbol{\theta}_h)\).

Throughout this section, the hypernetwork \(\mathcal H\) is implemented as a
one-dimensional CNN with six convolutional layers, 32 channels per layer, and
kernel size \(5\), together with boundary-consistent padding. In the
unknown-flux setting, \textsf{FluxNet} is implemented as a one-dimensional CNN
with four convolutional layers, 32 channels per layer, and kernel size \(5\),
again using padding matched to the boundary condition. For the Euler example,
\textsf{FluxNet} is taken to be pointwise, with \(1\times1\) kernels in every
layer.

\subsection{Inviscid Burgers equation}

We consider the scalar inviscid Burgers equation
\begin{equation}\label{eq:burgers}
u_t+\Big(\tfrac12 u^2\Big)_x=0,\qquad x\in(0,2\pi),\; t>0,
\end{equation}
with periodic boundary conditions. We study two families of initial data: a
smooth family that steepens into a single shock and a piecewise-constant family
with two initial discontinuities.

\subsubsection{Single-shock family}

The single-shock family is defined by
\begin{equation}\label{eq:burgers-singleshock-ic}
u(x,0)=a+b\sin x,
\qquad
a\in[-\epsilon,\epsilon],\quad
b\in[1-\epsilon,1+\epsilon],
\qquad
\epsilon=0.25.
\end{equation}
These initial data evolve from a smooth profile into a single shock.

\paragraph{Training and reference setup}
For this family, the training set consists of 200 trajectories generated by
sampling
\[
a_i\sim\mathcal U[-\epsilon,\epsilon],\qquad
b_i\sim\mathcal U[1-\epsilon,1+\epsilon],\qquad i=1,\dots,200,
\]
and evolving each realization to \(T=1.5\) with CFL\(=0.4\). Training is carried
out on the mesh family \(N=32l\), \(l=1,\dots,8\), that is,
\(N\in\{32,64,\dots,256\}\), with window length \(L=20\) and unroll depth
\(K=4\). Reference solutions are computed by classical WENO5 on \(N=512\). Unless stated
otherwise, the accuracy, conservation, and refinement results reported below use
the fixed test parameters
\(
(a,b)=(-0.062730,\,0.965973)
\)
in \eqref{eq:burgers-singleshock-ic}.

\begin{figure}[htbp]
  \centering
  \begin{subfigure}{0.4\textwidth}
    \centering
    \includegraphics[width=\textwidth]{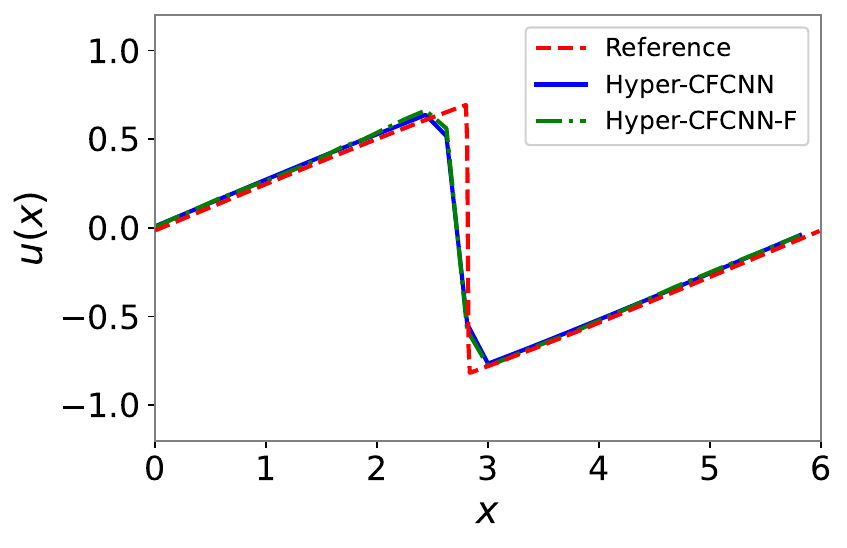}
    \caption{$N=32,t=3$}
    \label{fig:T3_N32}
  \end{subfigure}
  \hfill
  \begin{subfigure}{0.4\textwidth}
    \centering
    \includegraphics[width=\textwidth]{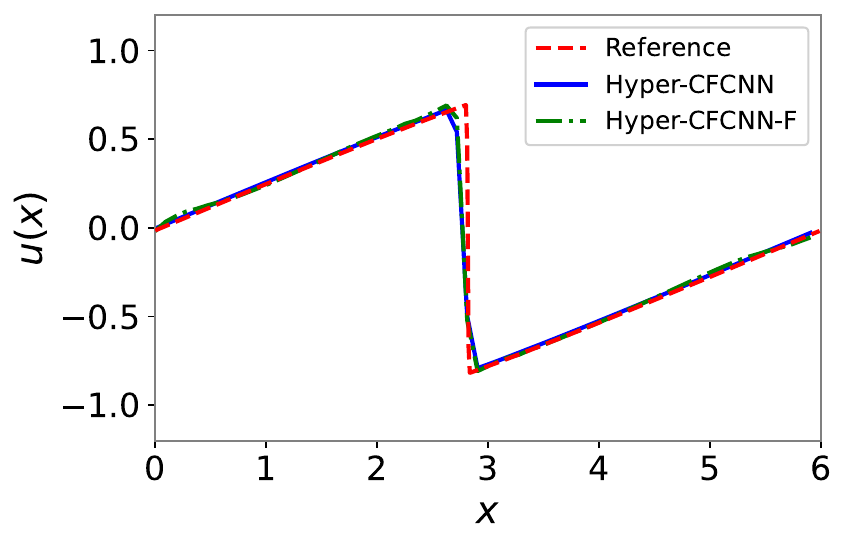}
    \caption{$N=64,t=3$}
    \label{fig:T3_N64}
  \end{subfigure}
  \vspace{1ex}
  \begin{subfigure}{0.4\textwidth}
    \centering
    \includegraphics[width=\textwidth]{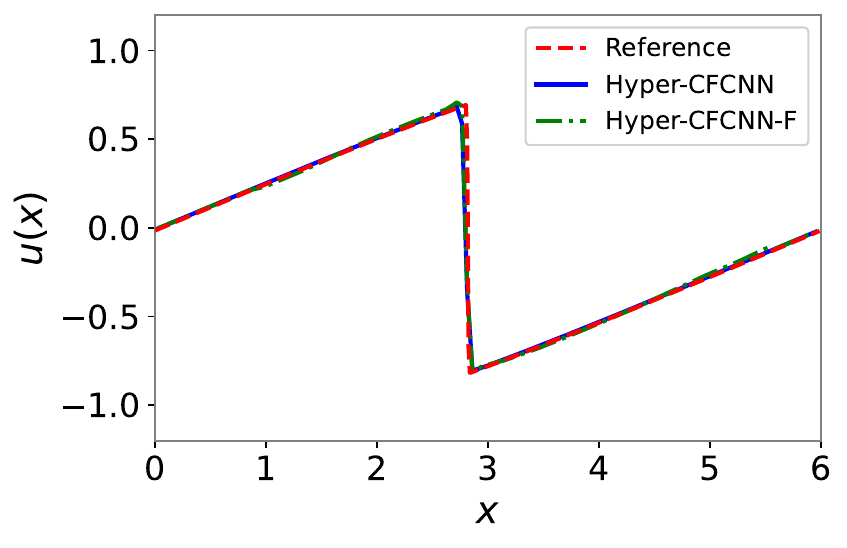}
    \caption{$N=128, t=3$}
    \label{fig:T3_N128}
  \end{subfigure}
  \hfill
  \begin{subfigure}{0.4\textwidth}
    \centering
    \includegraphics[width=\textwidth]{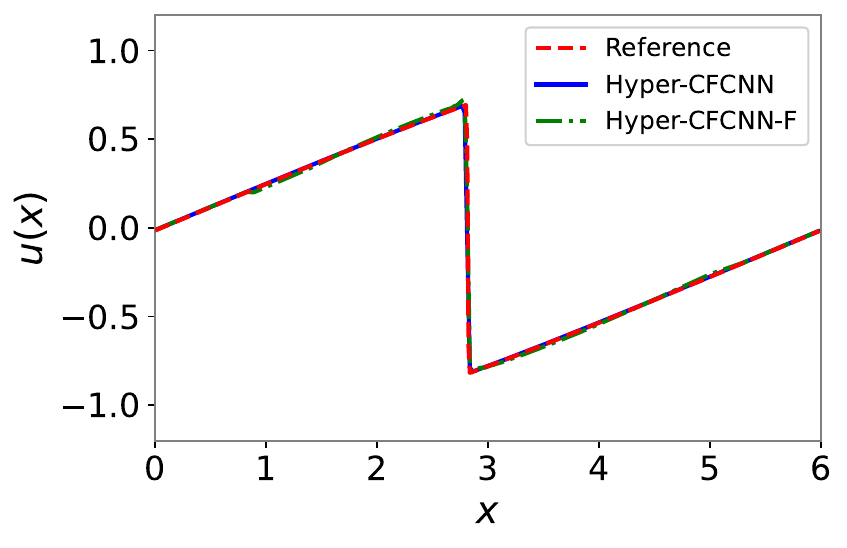}
    \caption{$N=256,t=3$}
    \label{fig:T3_N256}
  \end{subfigure}
  \caption{H--CFCNN and H--CFCNN--F versus the WENO5 reference at \(T=3\) for
the single-shock Burgers example on \(N=32,64,128,256\).}
  \label{fig:T3_mesh}
\end{figure}

\begin{figure}[htbp]
  \centering
  \includegraphics[width=0.6\textwidth]{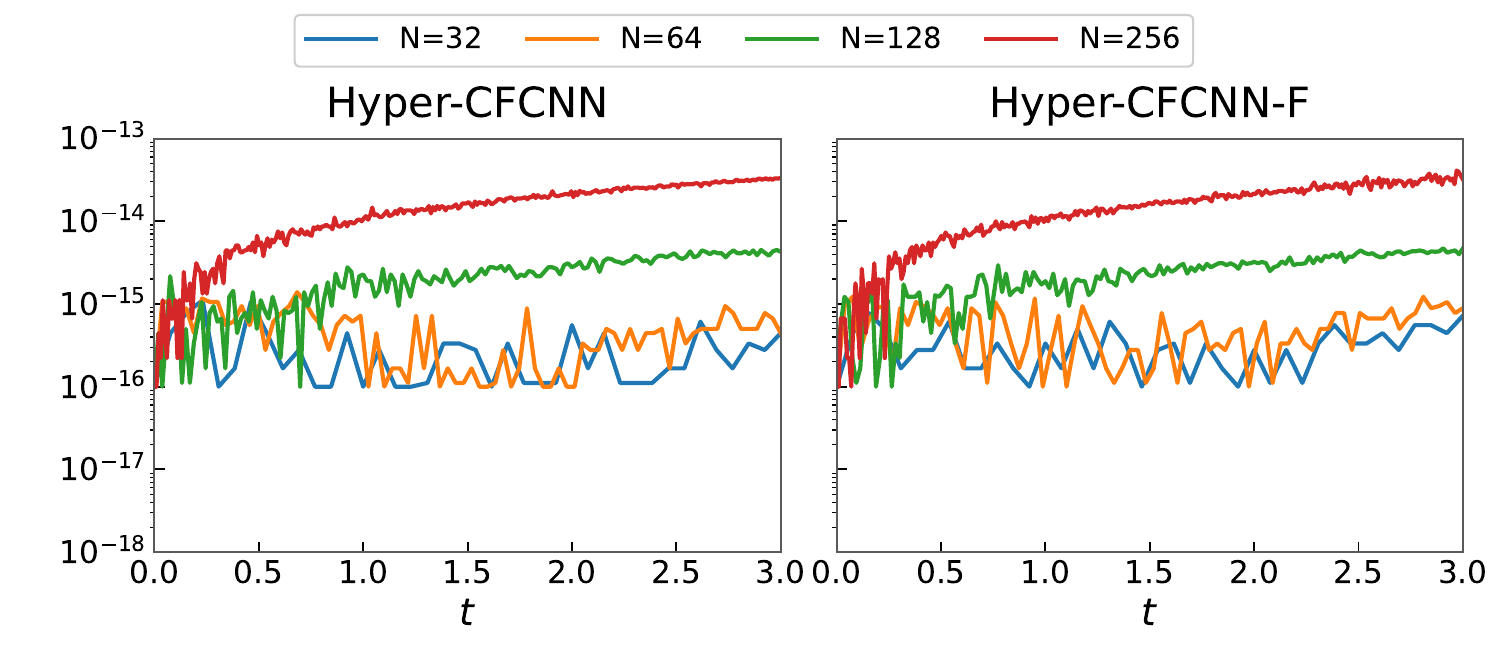}
 \caption{Conservation remainder \(C(u)\) versus time for H--CFCNN and
H--CFCNN--F on \(N=32,64,128,256\) for the single-shock Burgers example.}
  \label{fig:bur-mass}
\end{figure}

\paragraph{Accuracy, conservation, and refinement}
{\sloppy Figure~\ref{fig:T3_mesh} compares H--CFCNN, H--CFCNN--F, and the WENO5
reference at \(T=3\) for \(N\in\{32,64,128,256\}\). Both learned solvers improve
systematically under mesh refinement. On the finer meshes, the shock location
and surrounding profile are close to the reference solution, and the phase error
is smaller than on the coarser meshes. Figure~\ref{fig:bur-mass} reports the
conservation remainder \(C(u)\) up to \(T=3\) for \(N=32,64,128,256\). In this
periodic setting, both variants remain near machine precision throughout the
rollout, which is consistent with the conservative flux-difference form of the
update. \par}

Table~\ref{tab:error_T15} reports the mean-squared error (MSE) at \(T=1.5\) for
H--CFCNN together with the empirical refinement rate
\(p=\log_2(E_h/E_{h/2})\). The error decreases monotonically as the mesh is
refined, and the observed rate increases on the finer meshes. Table
\ref{tab:parameter_cost} reports the corresponding target-network parameter
count and wall-clock prediction time from \(t=0\) to \(t=3\). Both quantities
grow linearly with the number of cells, consistent with the per-cell parameter
generation induced by the hypernetwork.

\begin{figure}[htbp]
  \centering
  \begin{subfigure}{0.22\textwidth}
    \includegraphics[width=\textwidth]{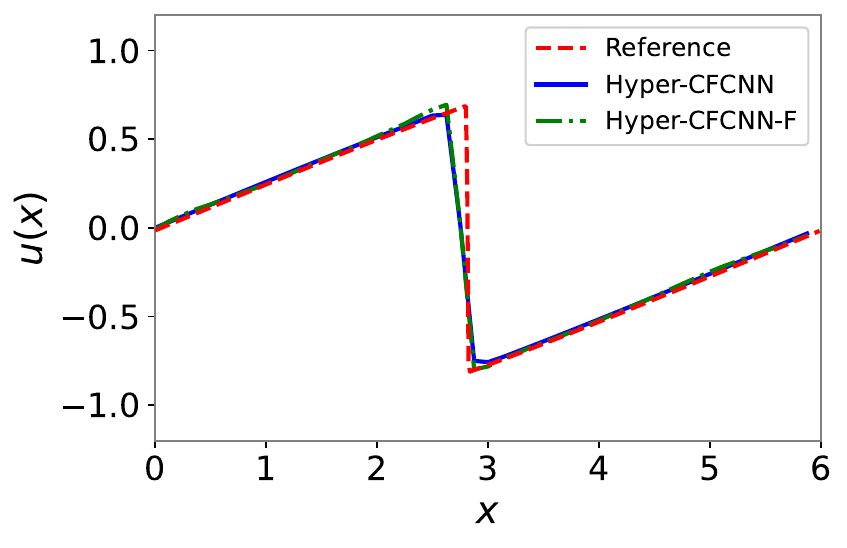}
    \caption{$N=48,T=3$ }
    \label{fig:T3_N48}
  \end{subfigure}
  \hfill
  \begin{subfigure}{0.22\textwidth}
    \includegraphics[width=\textwidth]{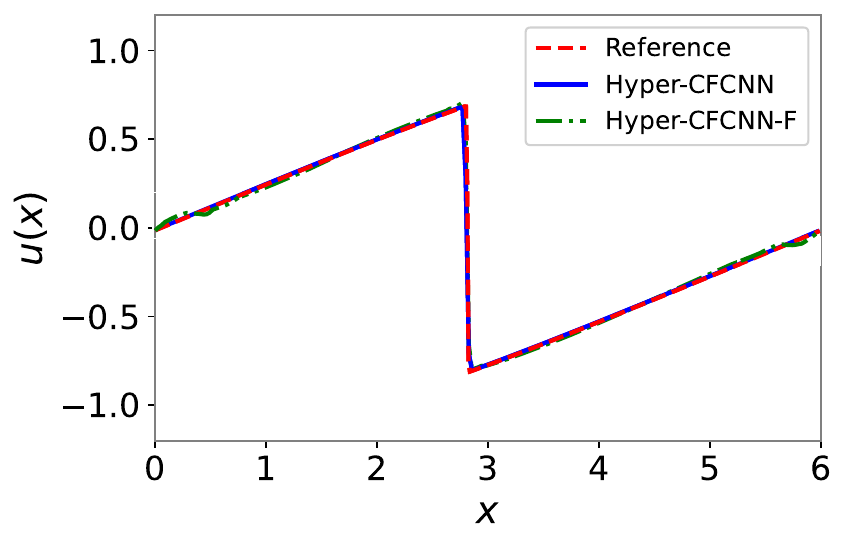}
    \caption{$N=208,T=3$}
    \label{fig:T3_N208}
  \end{subfigure}
  \hfill
  \begin{subfigure}{0.22\textwidth}
    \includegraphics[width=\textwidth]{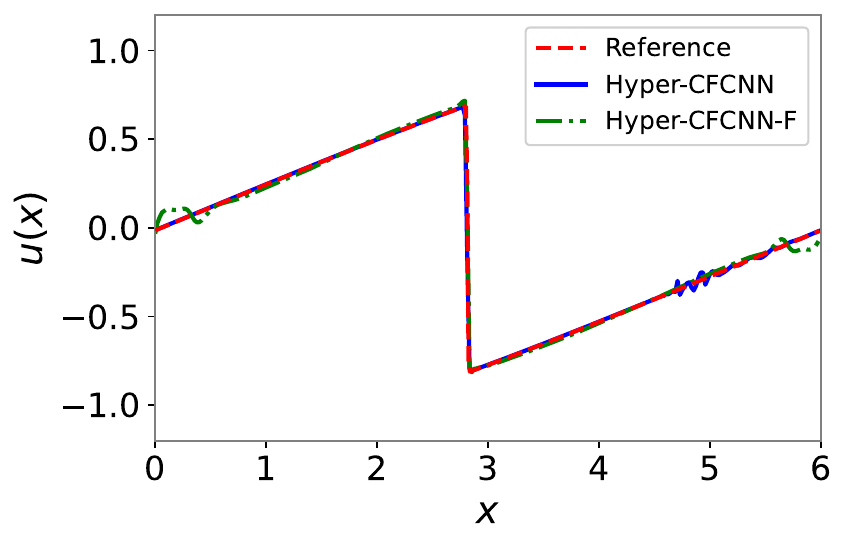}
    \caption{$N=288,T=3$}
    \label{fig:T3_N288}
  \end{subfigure}
  \hfill
  \begin{subfigure}{0.22\textwidth}
    \includegraphics[width=\textwidth]{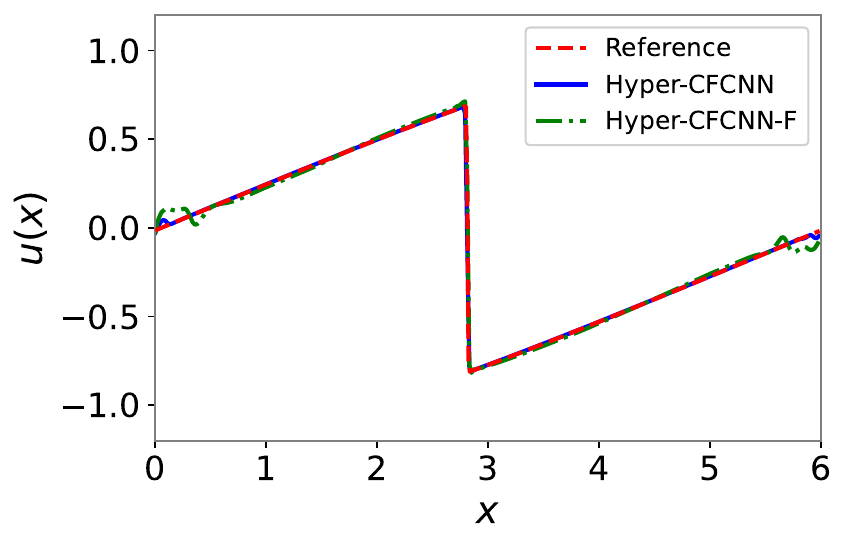}
    \caption{$N=320,T=3$ }
    \label{fig:T3_N320}
  \end{subfigure}

  \caption{H--CFCNN and H--CFCNN--F versus the reference at \(T=3\) for the
single-shock Burgers example on \(N=48,208,288,320\). The first two meshes lie
inside the training interval, whereas the last two lie outside it.}
  \label{fig:T3_extrap}
\end{figure}

\paragraph{Generalization across meshes and initial conditions}
Figure~\ref{fig:T3_extrap} examines transfer to unseen meshes both inside and
outside the trained resolution interval. The meshes \(N=48\) and \(N=208\) lie
inside the training range but are not used during training, whereas \(N=288\)
and \(N=320\) lie outside it. Both learned solvers remain close to the
reference on all four meshes, although the finest extrapolated cases show mild
degradation and a small phase lag.

\begin{figure}[htbp]
  \centering
  \begin{subfigure}{0.4\textwidth}
    \centering
    \includegraphics[width=\textwidth]{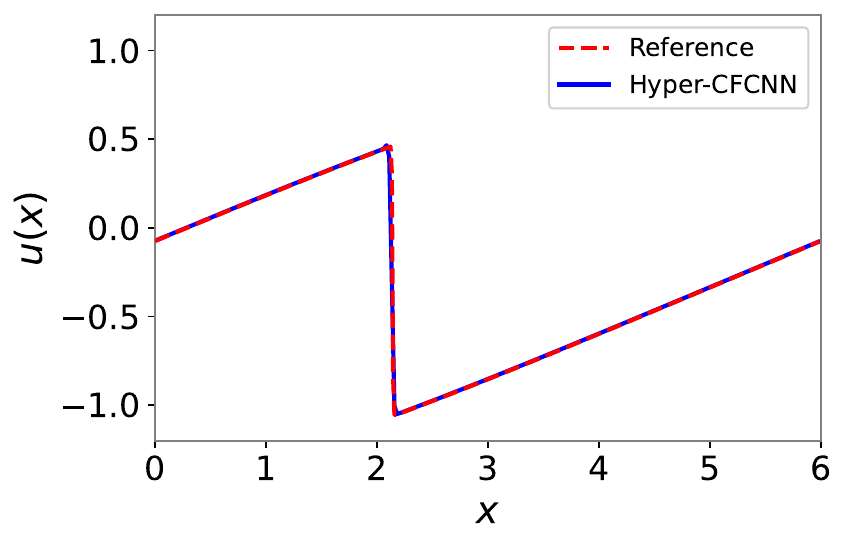}
    \caption{$a=-0.3,\ b=1.0$}
    \label{fig:init_a-03_b10}
  \end{subfigure}
  \hfill
  \begin{subfigure}{0.4\textwidth}
    \centering
    \includegraphics[width=\textwidth]{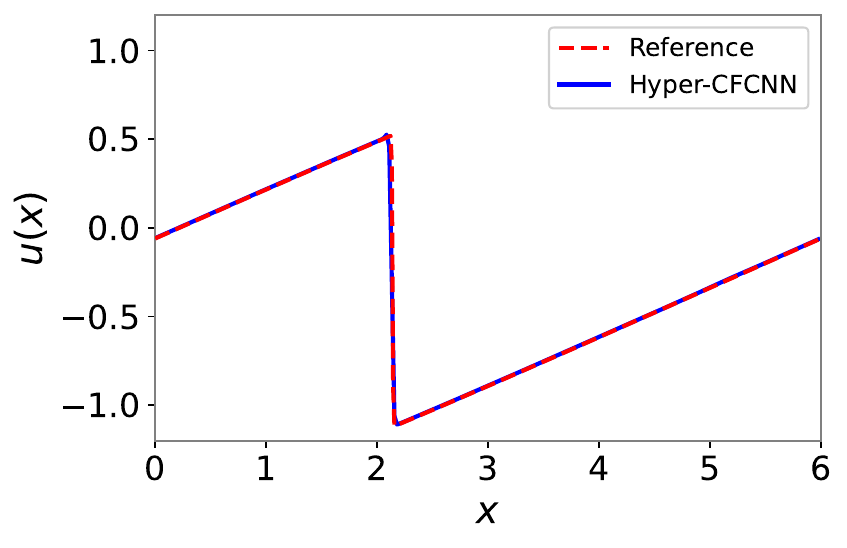}
    \caption{$a=-0.3,\ b=1.3$}
    \label{fig:init_a-03_b13}
  \end{subfigure}
 \caption{H--CFCNN versus the reference at \(T=3\) for two out-of-range initial
conditions on \(N=256\) for the single-shock Burgers example.}
  \label{fig:init_combined}
\end{figure}

Figure~\ref{fig:init_combined} shows two out-of-range initial conditions on
\(N=256\), corresponding to \((a,b)=(-0.3,1.0)\) and \((a,b)=(-0.3,1.3)\). In
the first case, \(a\) lies outside the training range, whereas in the second
case both \(a\) and \(b\) lie outside it. In both cases, H--CFCNN remains close
to the reference profile, which indicates good extrapolation to initial
conditions outside the training range.

\begin{figure}[htbp]
\centering
\begin{subfigure}{0.28\textwidth}
\includegraphics[width=\textwidth]{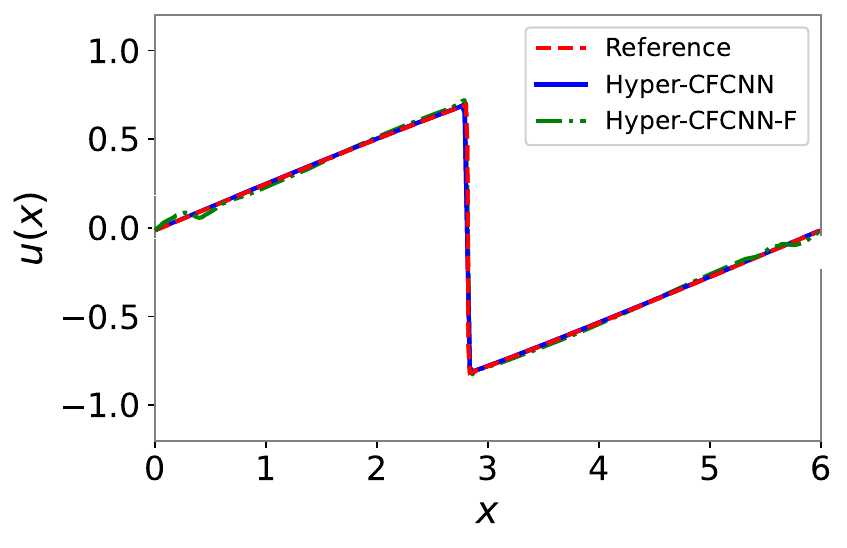}
\caption{$\Delta t = 0.2 \Delta x$}
\label{fig:dt_02}
\end{subfigure}
\hfill
\begin{subfigure}{0.28\textwidth}
\includegraphics[width=\textwidth]{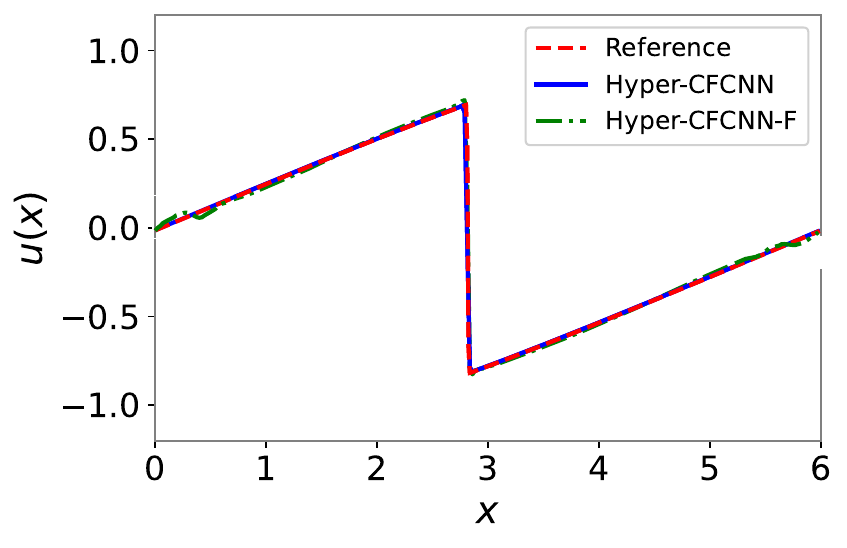}
\caption{$\Delta t = 0.3 \Delta x$}
\label{fig:dt_03}
\end{subfigure}
\hfill
\begin{subfigure}{0.28\textwidth}
\includegraphics[width=\textwidth]{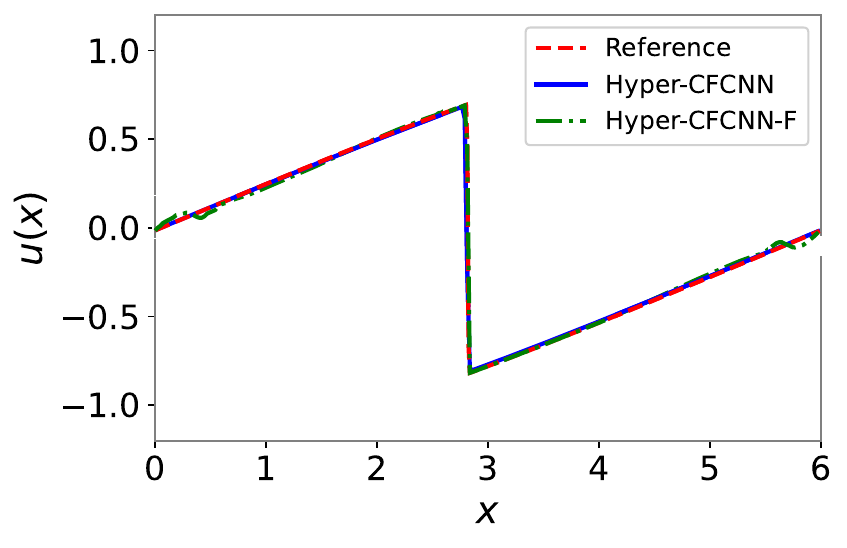}
\caption{$\Delta t = 0.48 \Delta x$}
\label{fig:dt_048}
\end{subfigure}
\caption{H--CFCNN and H--CFCNN--F versus the reference at \(T=3\) for
different time-step sizes on the \(N=256\) mesh for the single-shock Burgers
example.}
\label{fig:dt_flexibility}
\end{figure}

\paragraph{Sensitivity to the time step}
Figure~\ref{fig:dt_flexibility} reports the solution at \(T=3\) on \(N=256\)
for \(\Delta t/\Delta x\in\{0.2,0.3,0.48\}\), whereas training uses CFL\(=0.4\).
H--CFCNN remains accurate across these choices. H--CFCNN--F is also stable,
although small oscillations appear when the test step size differs more
substantially from the training value.

\begin{table}[!htbp]
\centering
\caption{MSE and empirical refinement rates for H--CFCNN on the single-shock Burgers example at \(T=1.5\).}
\label{tab:error_T15}
\begin{tabular}{ccc}
\toprule
Mesh size \(N\) & MSE & Order \(p\) \\
\midrule
32  & \(1.2034\times10^{-2}\) & --   \\
64  & \(5.5819\times10^{-3}\) & 0.55 \\
128 & \(1.5290\times10^{-3}\) & 0.93 \\
256 & \(2.9309\times10^{-4}\) & 1.19 \\
\bottomrule
\end{tabular}
\end{table}

\begin{table}[!htbp]
\centering
\caption{Target-network parameter count and wall-clock prediction time for H--CFCNN on the single-shock Burgers example from \(t=0\) to \(t=3\).}
\label{tab:parameter_cost}
\begin{tabular}{ccc}
\toprule
Mesh size \(N\) & Parameters & Time (s) \\
\midrule
32  & 2560  & 0.33 \\
64  & 5120  & 0.55 \\
128 & 10240 & 0.93 \\
256 & 20480 & 1.19 \\
\bottomrule
\end{tabular}
\end{table}

\subsubsection{Multi-shock family}

We again consider \eqref{eq:burgers} with periodic boundary conditions. The
initial condition is piecewise constant with two discontinuities,
\begin{equation}\label{eq:burgurs-shock}
u(x,0)=
\begin{cases}
y_1, & x\in[\min(x_1^0,x_2^0),\,\max(x_1^0,x_2^0)],\\[2pt]
y_2, & \text{otherwise},
\end{cases}
\end{equation}
where \(y_1,y_2\in[-1,1]\) and \(x_1^0,x_2^0\in[0,2\pi]\). 

\paragraph{Training and reference setup}
For this family, the training set consists of 200 trajectories generated by
sampling
\[
y_{i,1},y_{i,2}\sim\mathcal U[-1,1],\qquad
x^0_{i,1},x^0_{i,2}\sim\mathcal U[0,2\pi],\qquad i=1,\dots,200,
\]
and evolving each realization to \(T=1.5\) with CFL\(=0.4\). Training is carried
out on the mesh family \(N=32l\), \(l=1,\dots,8\), that is,
\(N\in\{32,64,\dots,256\}\), with window length \(L=20\) and unroll depth
\(K=4\). Reference solutions are again computed by classical WENO5 on \(N=512\). Unless stated otherwise, the accuracy, conservation, and long-time results reported
below are evaluated on the fixed test initial condition
\[
u(x,0)=
\begin{cases}
0.8,  & x\in[0,2.5)\cup[4.5,2\pi),\\
-0.1, & x\in[2.5,3.5),\\
-0.7, & x\in[3.5,4.5).
\end{cases}
\]
which lies outside the training family \eqref{eq:burgurs-shock} and produces a
more demanding shock-interaction pattern over a longer time interval.

\begin{figure}[htbp]
\centering

\begin{subfigure}{0.31\textwidth}
\includegraphics[width=\textwidth]{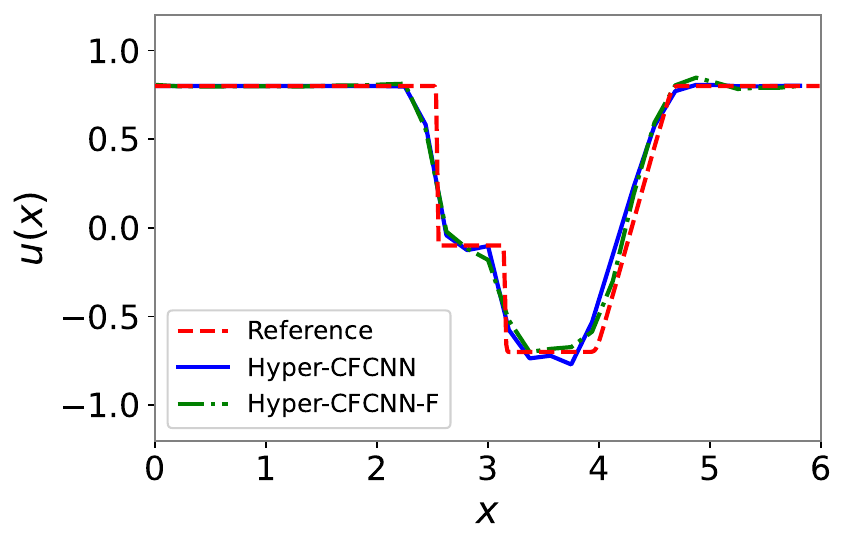}
\caption{$N=32,\ T=0.5$}
\label{fig:bur1-T05-N32}
\end{subfigure}
\hfill
\begin{subfigure}{0.31\textwidth}
\includegraphics[width=\textwidth]{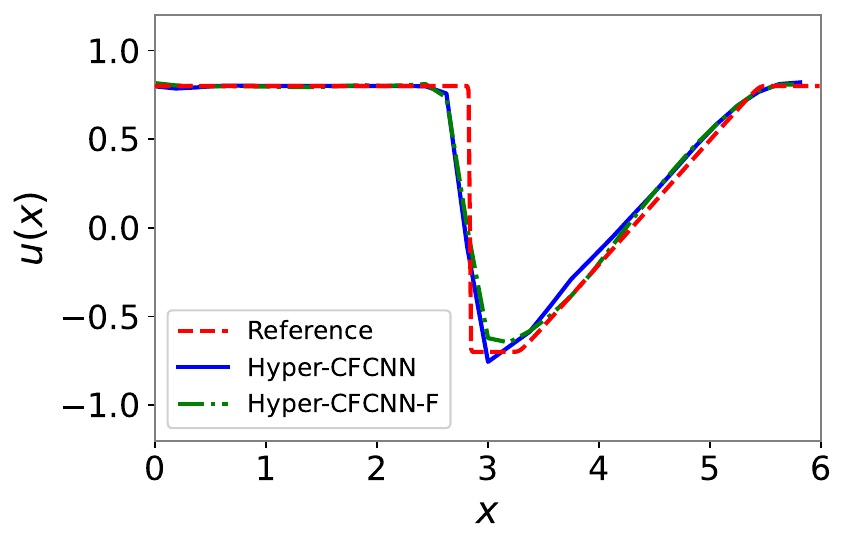}
\caption{$N=32,\ T=1.5$}
\label{fig:bur1-T15-N32}
\end{subfigure}
\hfill
\begin{subfigure}{0.31\textwidth}
\includegraphics[width=\textwidth]{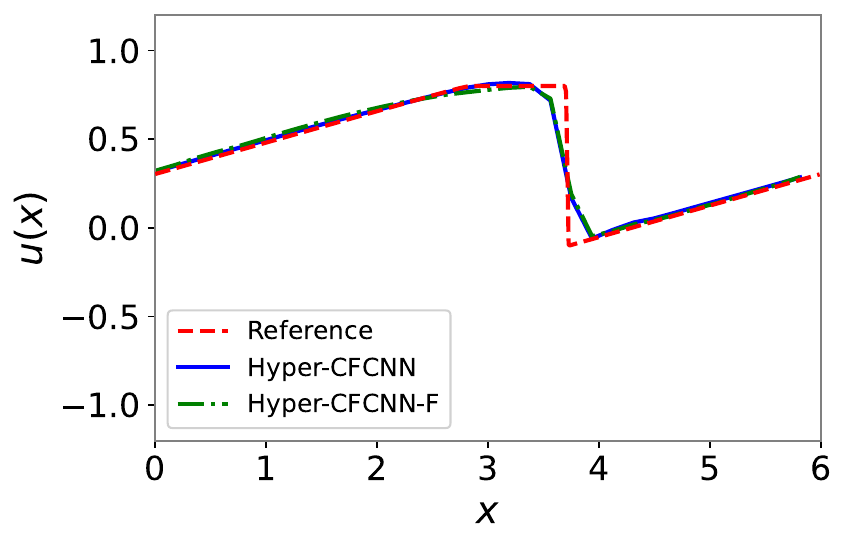}
\caption{$N=32,\ T=6$}
\label{fig:bur1-T6-N32}
\end{subfigure}

\vspace{0.8ex}

\begin{subfigure}{0.31\textwidth}
\includegraphics[width=\textwidth]{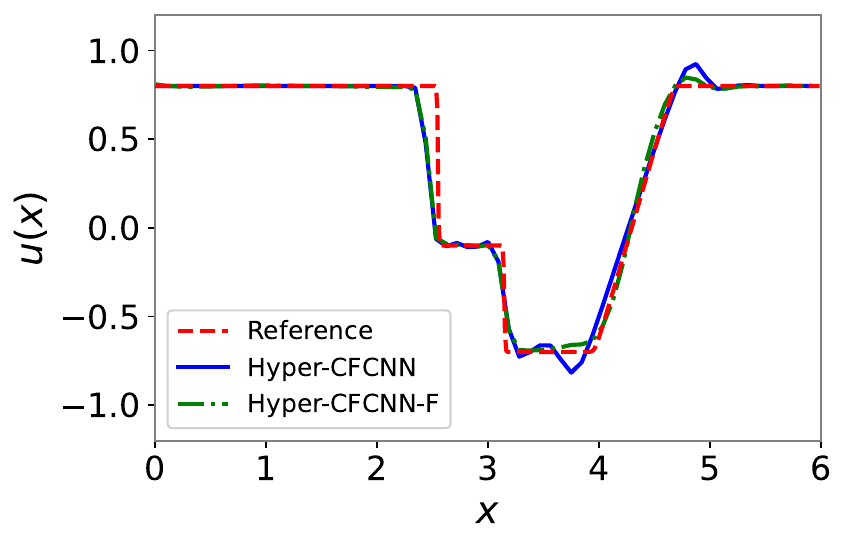}
\caption{$N=64,\ T=0.5$}
\label{fig:bur1-T05-N64}
\end{subfigure}
\hfill
\begin{subfigure}{0.31\textwidth}
\includegraphics[width=\textwidth]{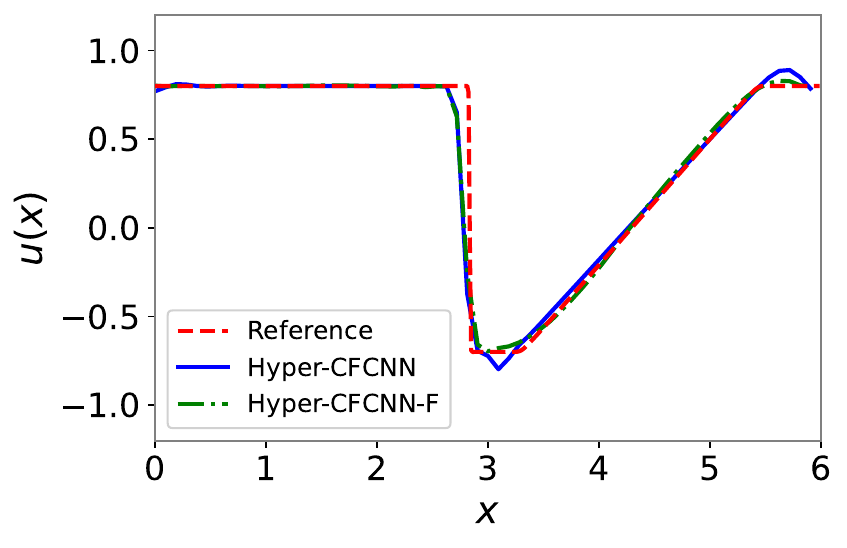}
\caption{$N=64,\ T=1.5$}
\label{fig:bur1-T15-N64}
\end{subfigure}
\hfill
\begin{subfigure}{0.31\textwidth}
\includegraphics[width=\textwidth]{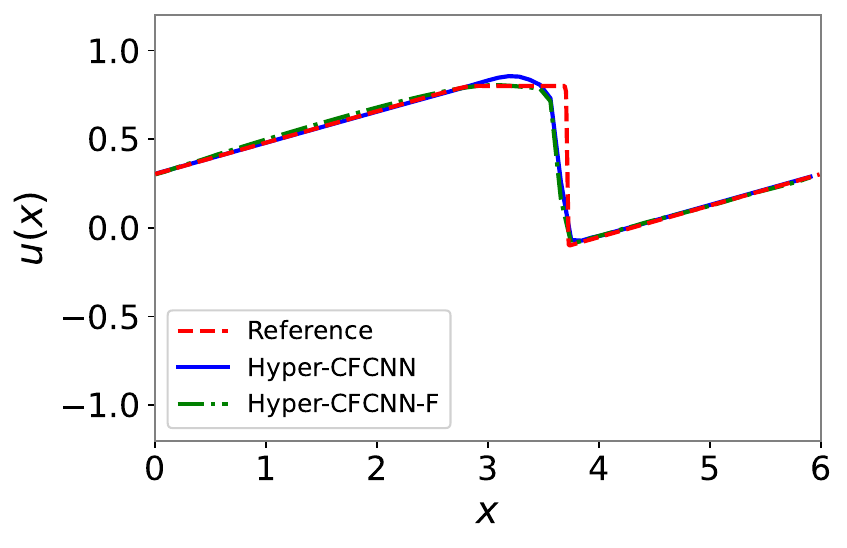}
\caption{$N=64,\ T=6$}
\label{fig:bur1-T6-N64}
\end{subfigure}

\vspace{0.8ex}

\begin{subfigure}{0.31\textwidth}
\includegraphics[width=\textwidth]{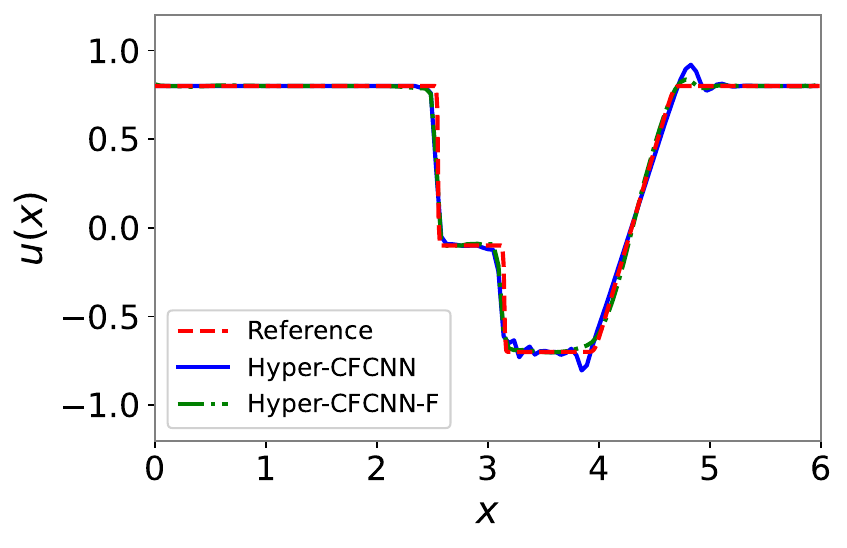}
\caption{$N=128,\ T=0.5$}
\label{fig:bur1-T05-N128}
\end{subfigure}
\hfill
\begin{subfigure}{0.31\textwidth}
\includegraphics[width=\textwidth]{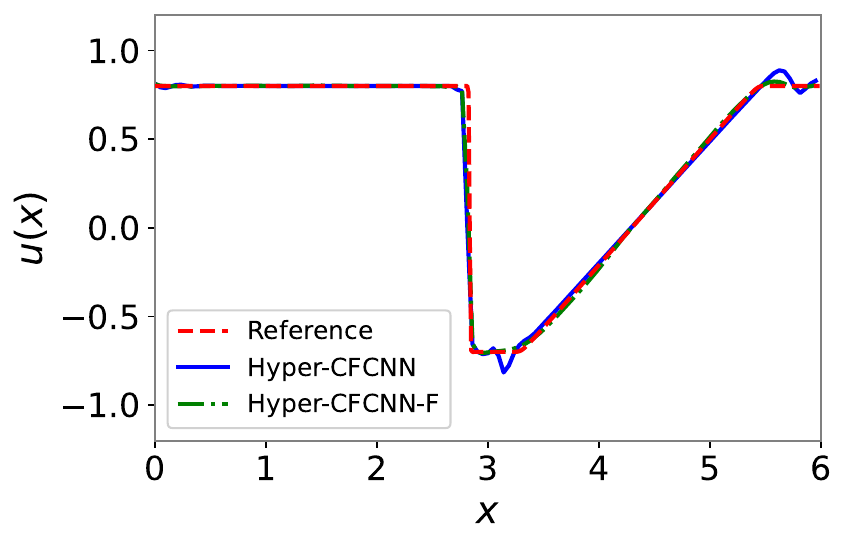}
\caption{$N=128,\ T=1.5$}
\label{fig:bur1-T15-N128}
\end{subfigure}
\hfill
\begin{subfigure}{0.31\textwidth}
\includegraphics[width=\textwidth]{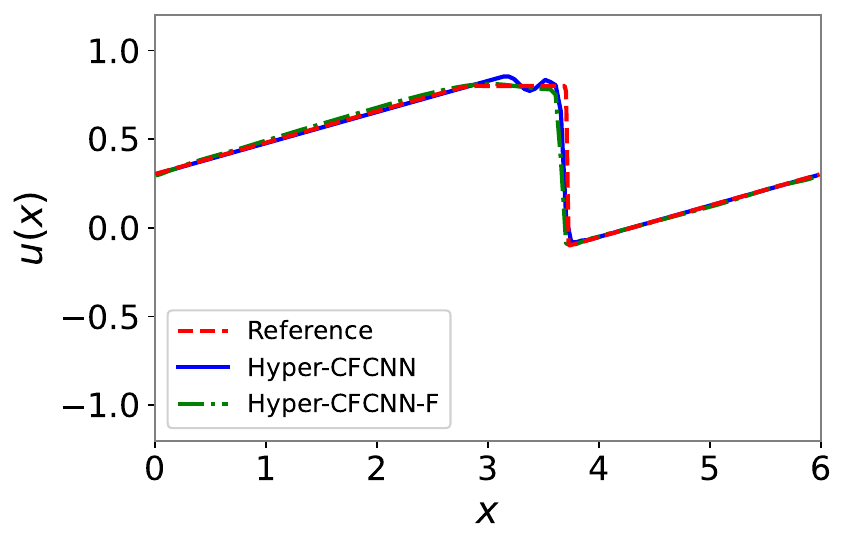}
\caption{$N=128,\ T=6$}
\label{fig:bur1-T6-N128}
\end{subfigure}

\vspace{0.8ex}

\begin{subfigure}{0.31\textwidth}
\includegraphics[width=\textwidth]{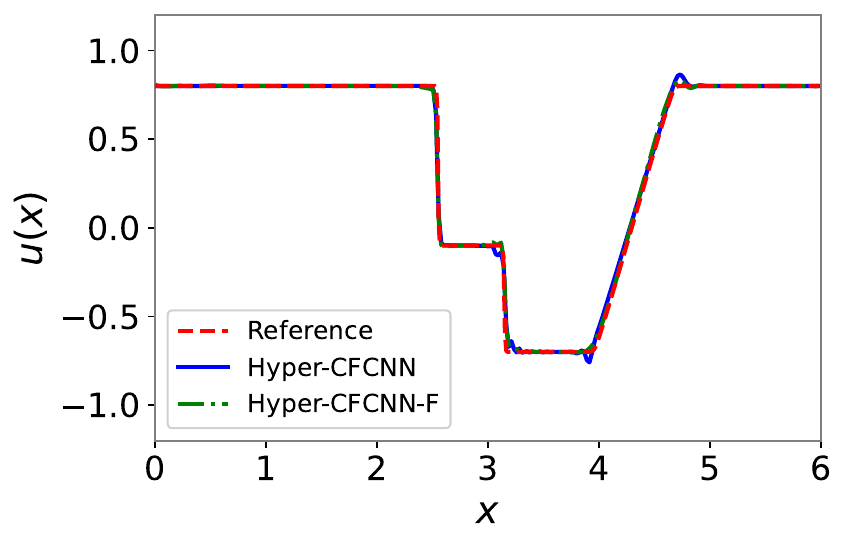}
\caption{$N=256,\ T=0.5$}
\label{fig:bur1-T05-N256}
\end{subfigure}
\hfill
\begin{subfigure}{0.31\textwidth}
\includegraphics[width=\textwidth]{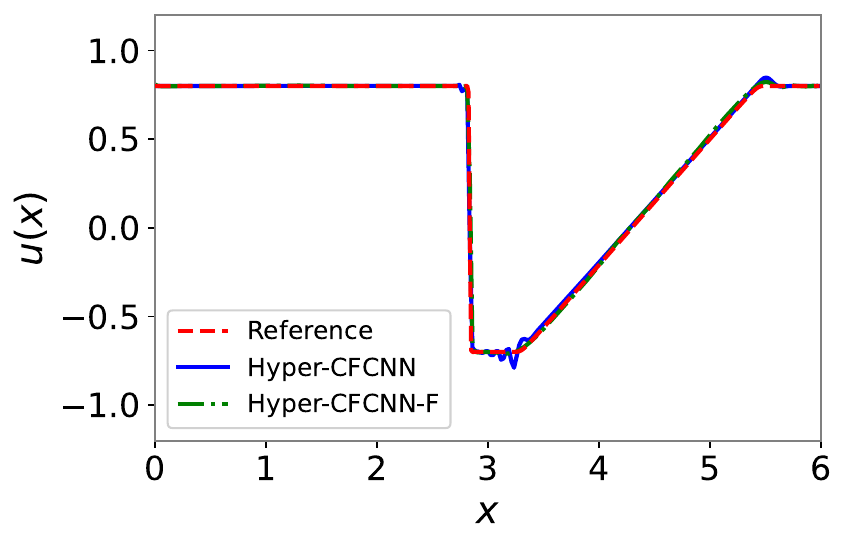}
\caption{$N=256,\ T=1.5$}
\label{fig:bur1-T15-N256}
\end{subfigure}
\hfill
\begin{subfigure}{0.31\textwidth}
\includegraphics[width=\textwidth]{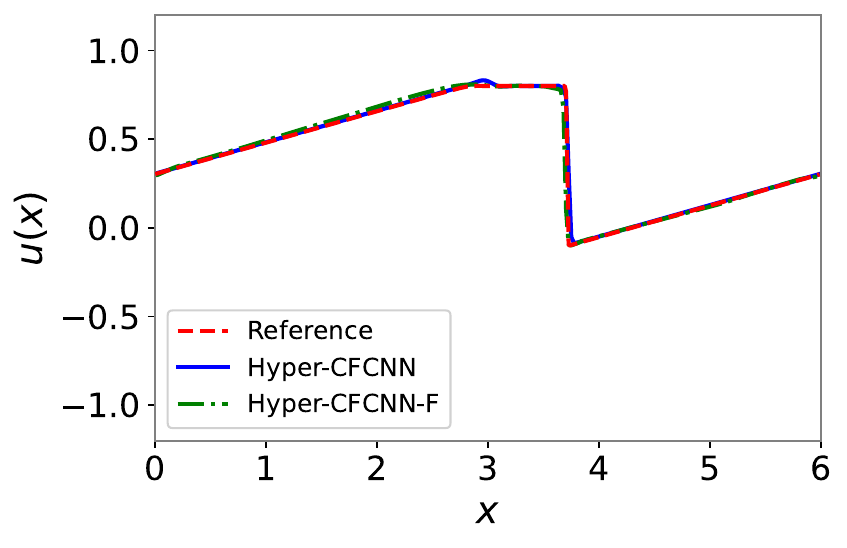}
\caption{$N=256,\ T=6$}
\label{fig:bur1-T6-N256}
\end{subfigure}

\caption{H--CFCNN and H--CFCNN--F versus the reference for the multi-shock
Burgers example at \(T=0.5\), \(1.5\), and \(6.0\) on \(N=32,64,128,256\).}
\label{fig:bur1-1}
\end{figure}

\begin{figure}[htbp]
  \centering
  \includegraphics[width=0.6\textwidth]{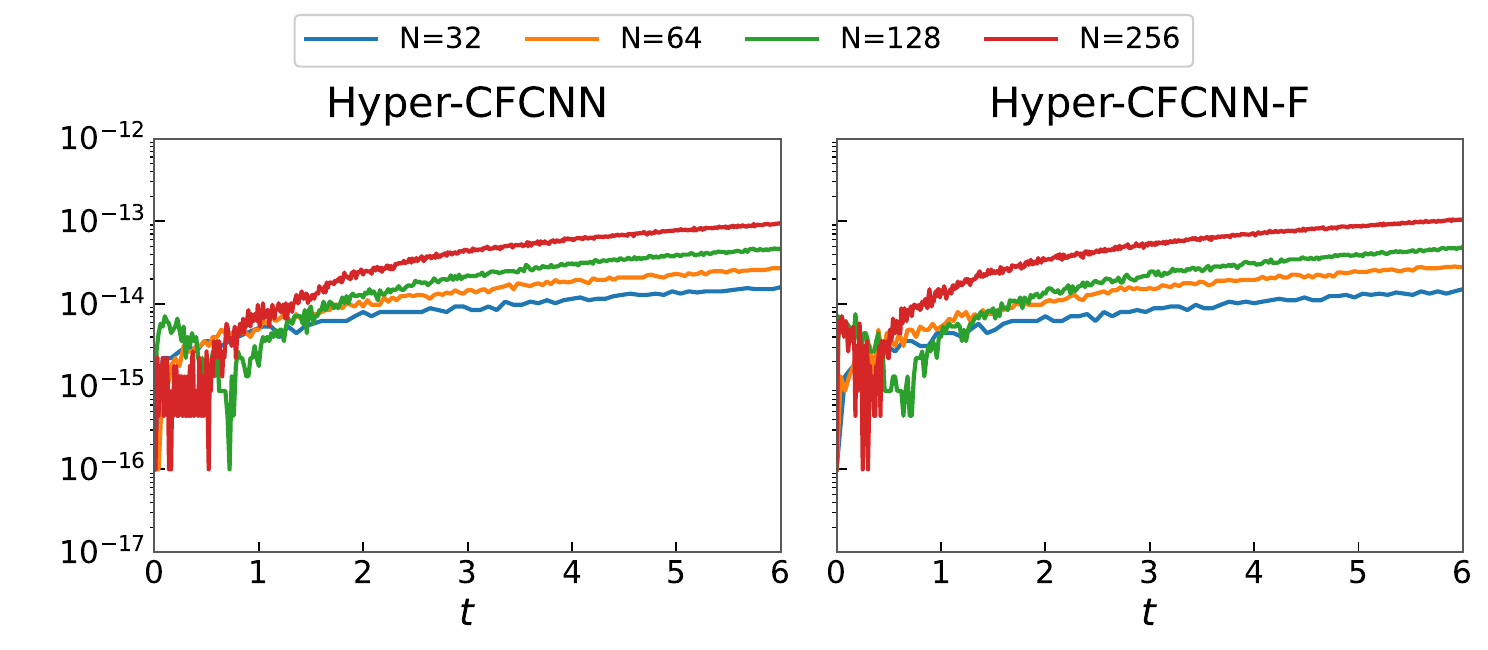}
  \caption{Conservation remainder \(C(u)\) versus time for H--CFCNN and
H--CFCNN--F on \(N=32,64,128,256\) for the multi-shock Burgers example.}
  \label{fig:bur1-mass}
\end{figure}

\paragraph{Accuracy, conservation, and long-time behavior}
Figure~\ref{fig:bur1-1} compares H--CFCNN and H--CFCNN--F with the reference at
\(T=0.5,1.5,6.0\) on \(N=32,64,128,256\). At \(T=0.5\), both methods sharply
resolve the propagating discontinuities. At \(T=1.5\), the interaction region is
well captured, and the discrepancy with the reference decreases under mesh
refinement. At \(T=6.0\), which lies well beyond the training horizon, both
variants still remain close to the reference and reproduce the post-collision
structure with good fidelity. Figure~\ref{fig:bur1-mass} shows the conservation remainder \(C(u)\) up to
\(T=6\) for \(N=32,64,128,256\). For both H--CFCNN and H--CFCNN--F, the
remainder remains near machine precision throughout the rollout.

\begin{figure}[htbp]
\centering
\begin{subfigure}{0.28\textwidth}
\includegraphics[width=\textwidth]{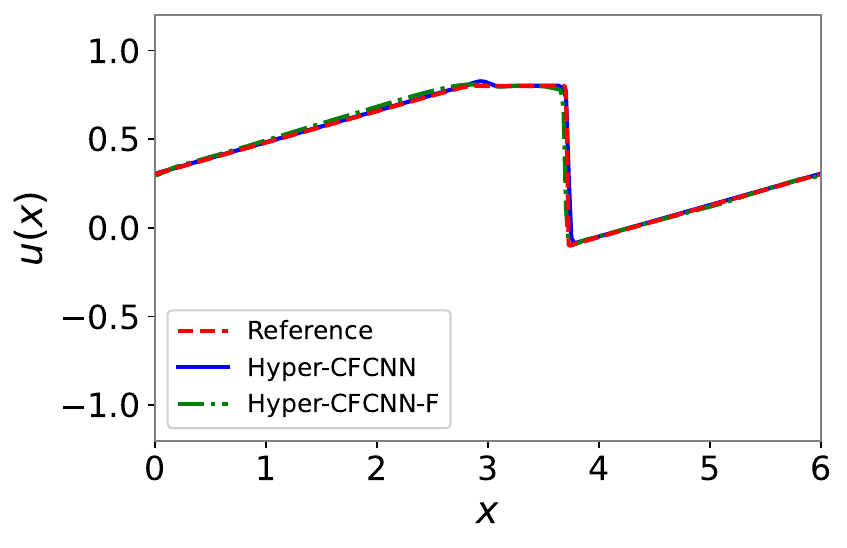}
\caption{$\Delta t = 0.2 \Delta x$}
\label{fig:bur1-dt_02}
\end{subfigure}
\hfill
\begin{subfigure}{0.28\textwidth}
\includegraphics[width=\textwidth]{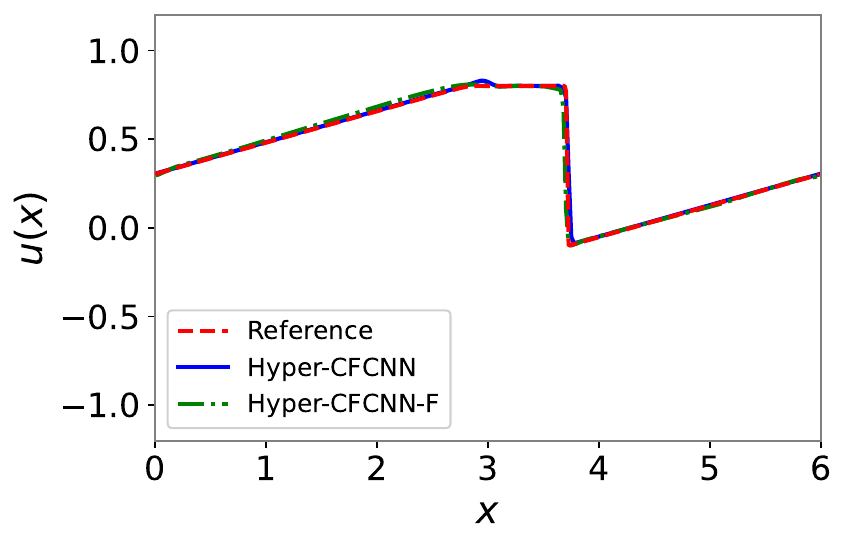}
\caption{$\Delta t = 0.3 \Delta x$}
\label{fig:bur1-dt_03}
\end{subfigure}
\hfill
\begin{subfigure}{0.28\textwidth}
\includegraphics[width=\textwidth]{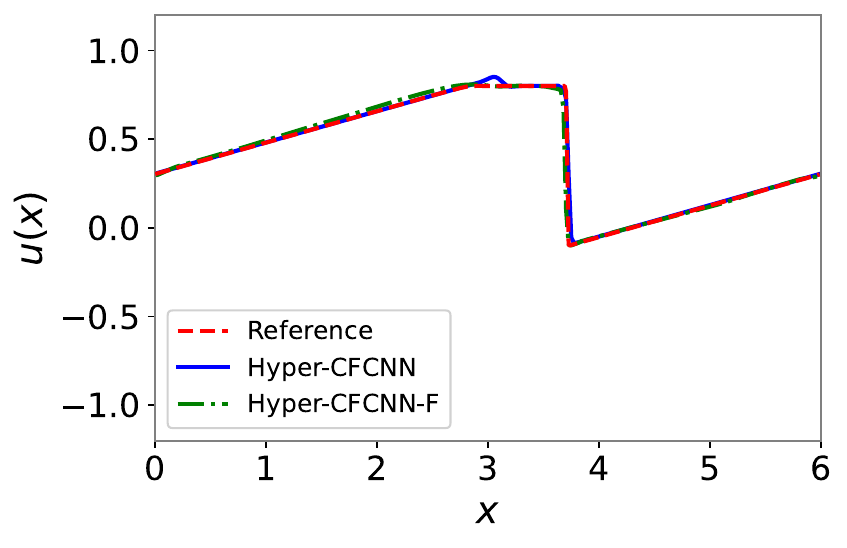}
\caption{$\Delta t = 0.48 \Delta x$}
\label{fig:bur1-dt_048}
\end{subfigure}
\caption{H--CFCNN and H--CFCNN--F versus the reference at \(T=6\) for
different time-step sizes on the \(N=256\) mesh for the multi-shock Burgers
example.}
\label{fig:bur1-dt_flexibility}
\end{figure}

\begin{figure}[htbp]
  \centering
  \begin{subfigure}{0.4\textwidth}
    \centering
    \includegraphics[width=\textwidth]{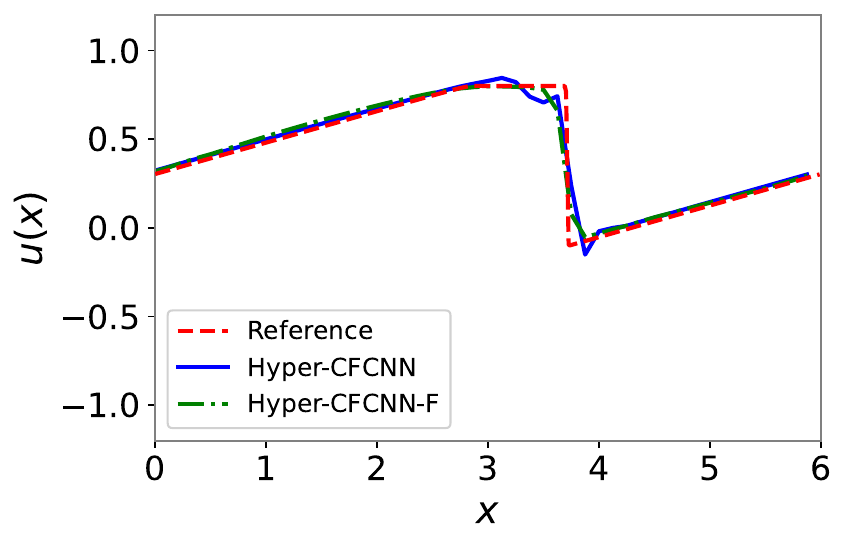}
    \caption{$N=48,\ T=6$}
    \label{fig:T6_N48}
  \end{subfigure}
  \hfill
  \begin{subfigure}{0.4\textwidth}
    \centering
    \includegraphics[width=\textwidth]{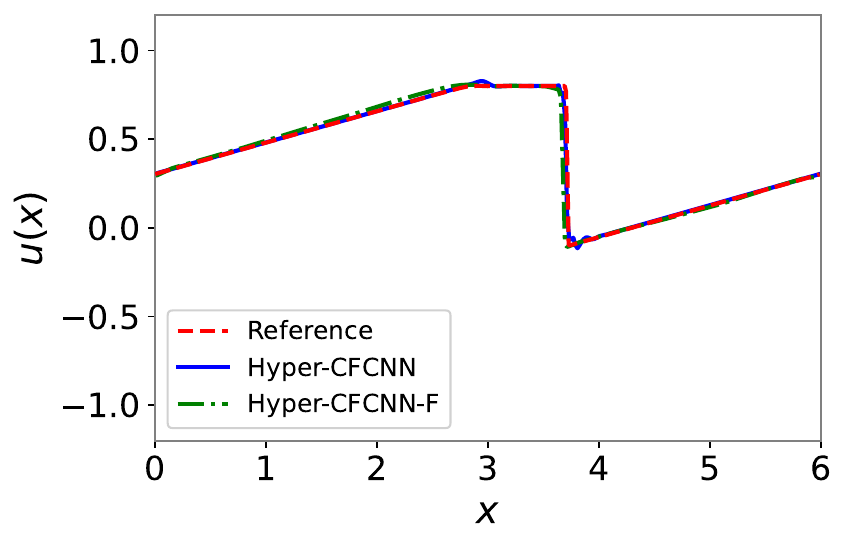}
    \caption{$N=320,\ T=6$}
    \label{fig:T6_N320}
  \end{subfigure}
  \caption{H--CFCNN and H--CFCNN--F versus the reference at \(T=6\) on the
unseen meshes \(N=48\) and \(N=320\) for the multi-shock Burgers example.}
  \label{fig:bur1_mesh_extrap}
\end{figure}

\paragraph{Time-step sensitivity and mesh transfer}
Figure~\ref{fig:bur1-dt_flexibility} reports the solution at \(T=6\) on
\(N=256\) for \(\Delta t/\Delta x\in\{0.2,0.3,0.48\}\). Both variants remain
stable over this range. Figure~\ref{fig:bur1_mesh_extrap} further shows
results at the unseen meshes \(N=48\) and \(N=320\). In both cases, the learned
solutions remain close to the reference, indicating good transfer to both
interpolatory and extrapolatory resolutions in this more demanding setting.

\subsection{Shallow-water equations}\label{sec:sw}

{\sloppy We consider the one-dimensional shallow-water system
\begin{equation}\label{eq:shallow}
\begin{aligned}
h_t + (hv)_x &= 0,\\
(hv)_t + \bigl(hv^2 + \tfrac{1}{2}gh^2\bigr)_x &= 0,
\end{aligned}
\end{equation}
on \(x\in(-5,5)\), with gravitational constant \(g=1\). We impose no-flux boundary conditions, implemented in the CNN components through
replicate padding as described in Section~\ref{sec:bc-fluxnet}. The initial
data are Riemann-type states,\par}
\begin{equation}\label{eq:shallow-initial}
h(x,0)=
\begin{cases}
h_\ell, & x\le x_0,\\
h_r,    & \text{otherwise},
\end{cases}
\qquad
v(x,0)=
\begin{cases}
v_\ell, & x\le x_0,\\
v_r,    & \text{otherwise},
\end{cases}
\end{equation}
where
\[
h_\ell\sim\mathcal U(2-\varepsilon_{h_\ell},\,2+\varepsilon_{h_\ell}),\qquad
h_r\sim\mathcal U(2-\varepsilon_{h_r},\,2+\varepsilon_{h_r}),
\]
and
\[
v_\ell,v_r,x_0\sim\mathcal U(-\varepsilon,\,\varepsilon),
\qquad
\varepsilon_{h_\ell}=0.2,\quad
\varepsilon_{h_r}=0.1,\quad
\varepsilon=0.1.
\]

\paragraph{Training and reference setup}
We generate 50 training trajectories by sampling the parameters
\((h_\ell,h_r,v_\ell,v_r,x_0)\) from the above ranges and evolving each case
from \(t=0\) to \(t=1\) with CFL\(=0.4\). Training is carried out on the mesh
family \(N\in\{64,128,192,256\}\), with window length \(L=30\) and unroll depth
\(K=4\). H--CFCNN uses the analytical shallow-water flux, whereas
H--CFCNN--F replaces it by \textsf{FluxNet}. Reference solutions are computed by WENO5 on \(N=1024\). Unless stated
otherwise, the accuracy, conservation, and refinement results reported below use
the fixed test parameters
\(
(h_\ell,h_r,v_\ell,v_r,x_0)=
(1.949816,\,
1.090143,\,
0.046399,\,
0.019732,\,
-0.068796)
\)
in \eqref{eq:shallow-initial}.

\begin{figure}[htbp]
\centering

\begin{subfigure}{0.31\textwidth}
\includegraphics[width=\textwidth]{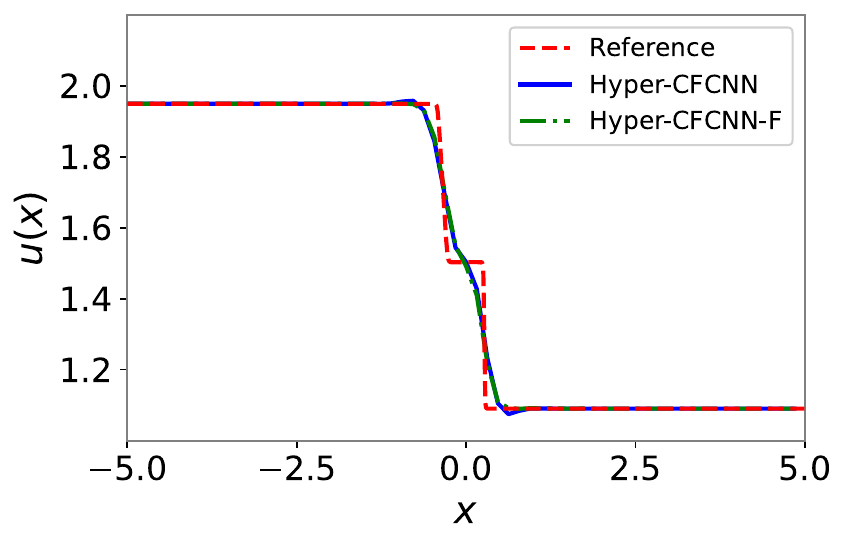}
\caption{$h$, $T=0.5$, $N=64$}
\label{fig:sw_h_T05_N64}
\end{subfigure}
\hfill
\begin{subfigure}{0.31\textwidth}
\includegraphics[width=\textwidth]{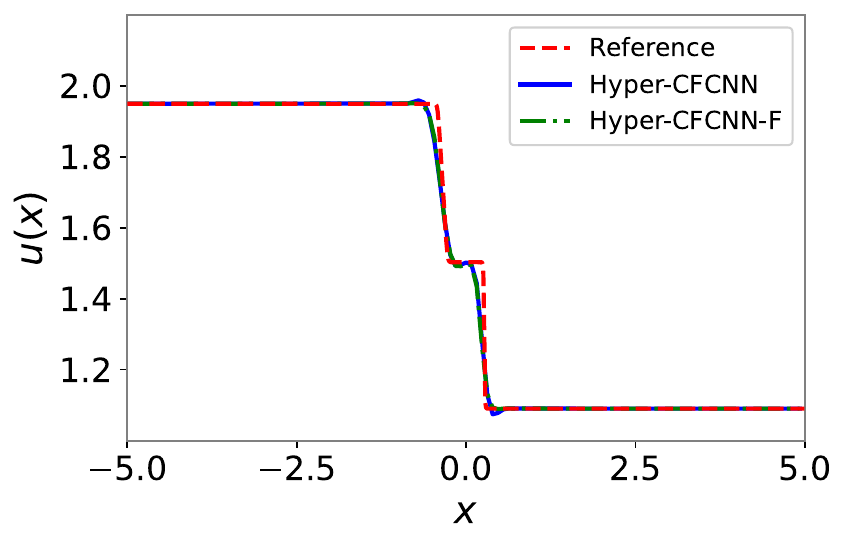}
\caption{$h$, $T=0.5$, $N=128$}
\label{fig:sw_h_T05_N128}
\end{subfigure}
\hfill
\begin{subfigure}{0.31\textwidth}
\includegraphics[width=\textwidth]{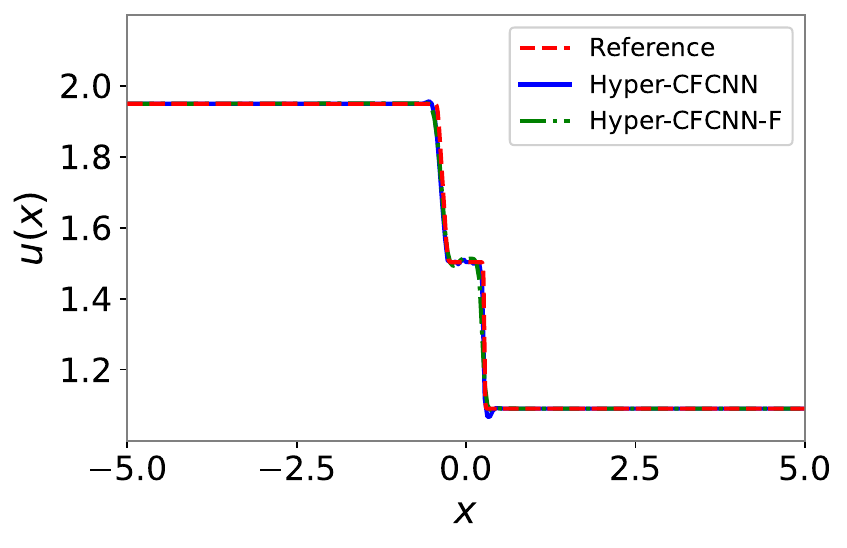}
\caption{$h$, $T=0.5$, $N=256$}
\label{fig:sw_h_T05_N256}
\end{subfigure}

\vspace{0.8ex}

\begin{subfigure}{0.31\textwidth}
\includegraphics[width=\textwidth]{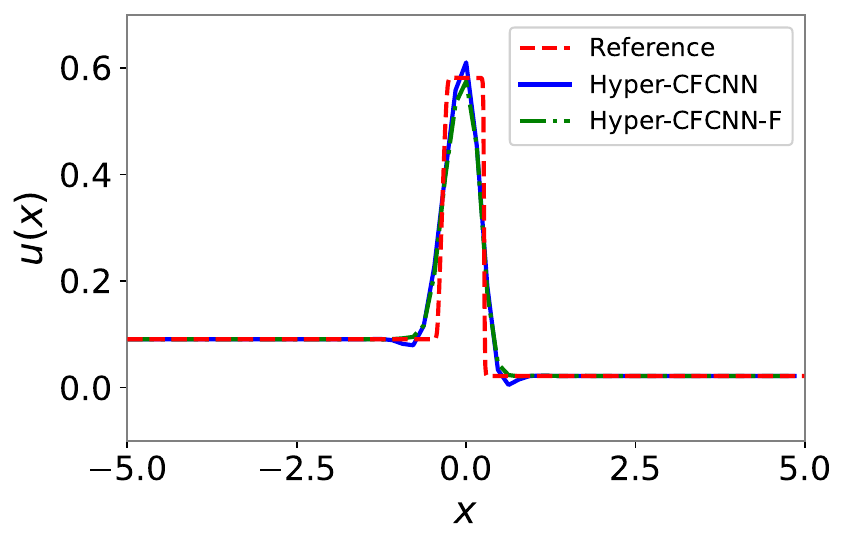}
\caption{$hv$, $T=0.5$, $N=64$}
\label{fig:sw_hv_T05_N64}
\end{subfigure}
\hfill
\begin{subfigure}{0.31\textwidth}
\includegraphics[width=\textwidth]{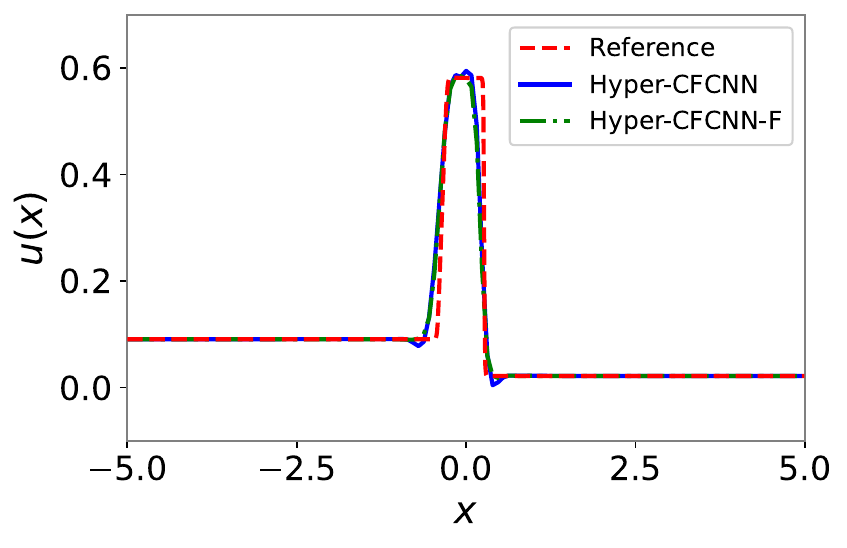}
\caption{$hv$, $T=0.5$, $N=128$}
\label{fig:sw_hv_T05_N128}
\end{subfigure}
\hfill
\begin{subfigure}{0.31\textwidth}
\includegraphics[width=\textwidth]{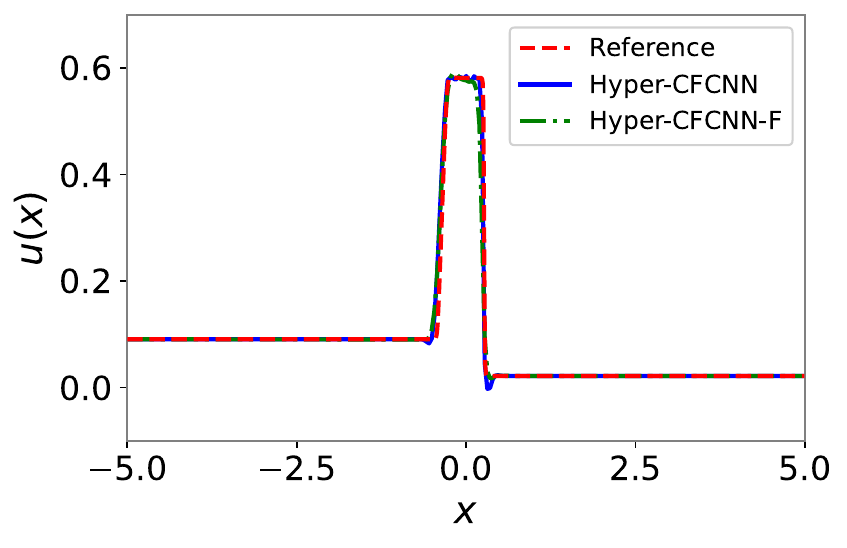}
\caption{$hv$, $T=0.5$, $N=256$}
\label{fig:sw_hv_T05_N256}
\end{subfigure}

\vspace{0.8ex}

\begin{subfigure}{0.31\textwidth}
\includegraphics[width=\textwidth]{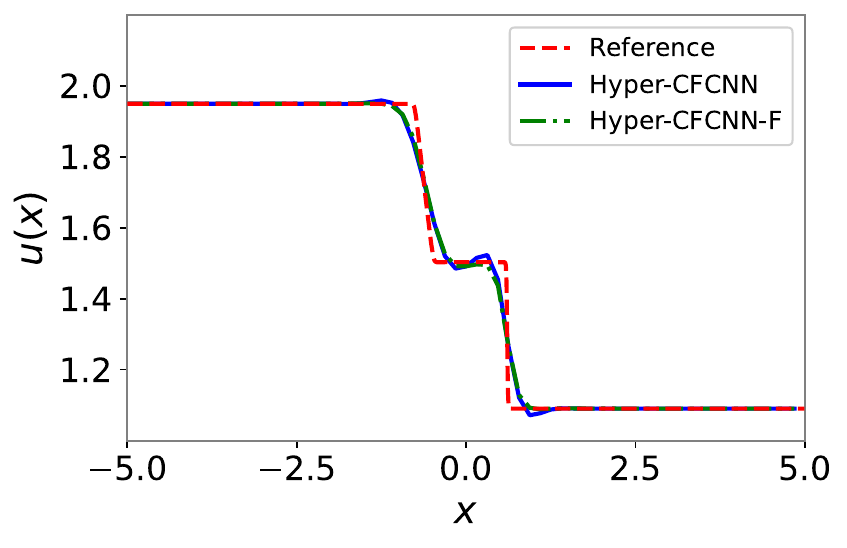}
\caption{$h$, $T=1.0$, $N=64$}
\label{fig:sw_h_T1_N64}
\end{subfigure}
\hfill
\begin{subfigure}{0.31\textwidth}
\includegraphics[width=\textwidth]{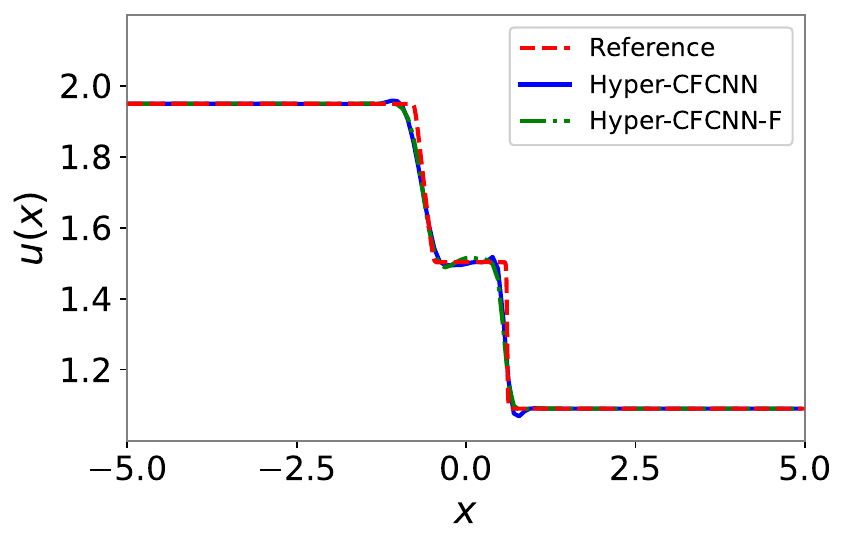}
\caption{$h$, $T=1.0$, $N=128$}
\label{fig:sw_h_T1_N128}
\end{subfigure}
\hfill
\begin{subfigure}{0.31\textwidth}
\includegraphics[width=\textwidth]{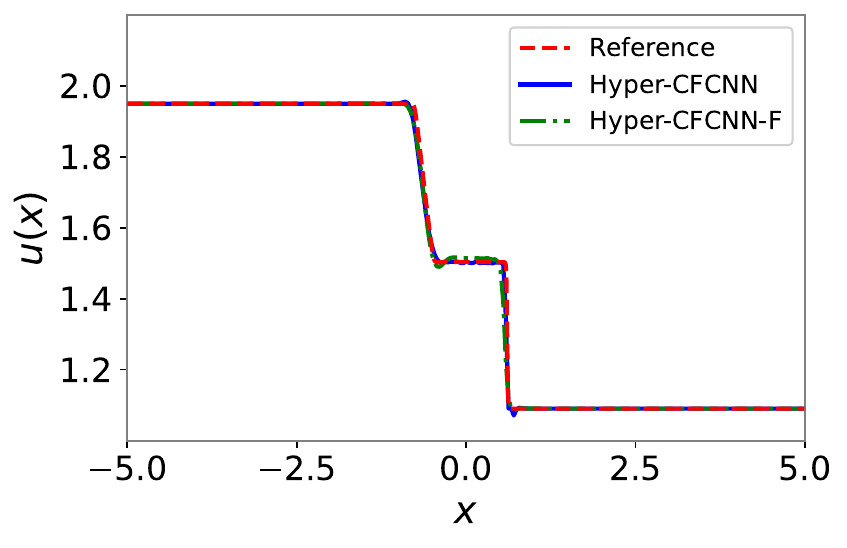}
\caption{$h$, $T=1.0$, $N=256$}
\label{fig:sw_h_T1_N256}
\end{subfigure}

\vspace{0.8ex}

\begin{subfigure}{0.31\textwidth}
\includegraphics[width=\textwidth]{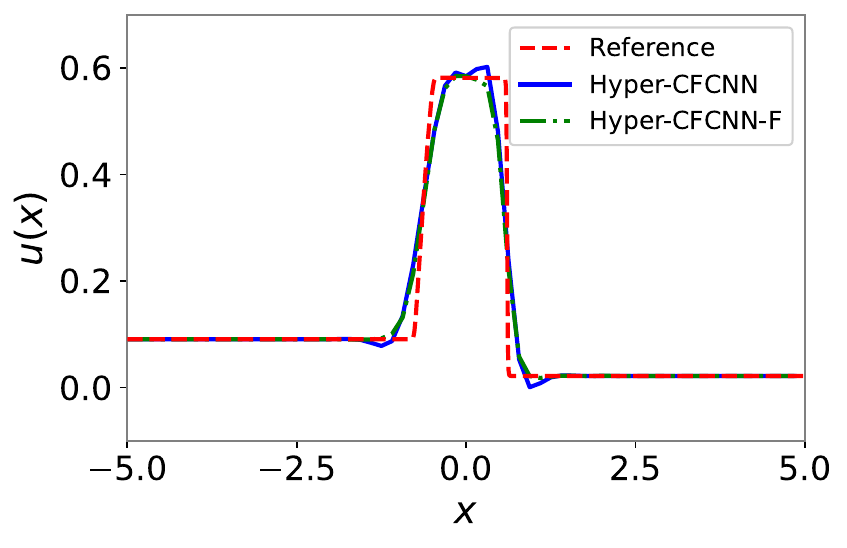}
\caption{$hv$, $T=1.0$, $N=64$}
\label{fig:sw_hv_T1_N64}
\end{subfigure}
\hfill
\begin{subfigure}{0.31\textwidth}
\includegraphics[width=\textwidth]{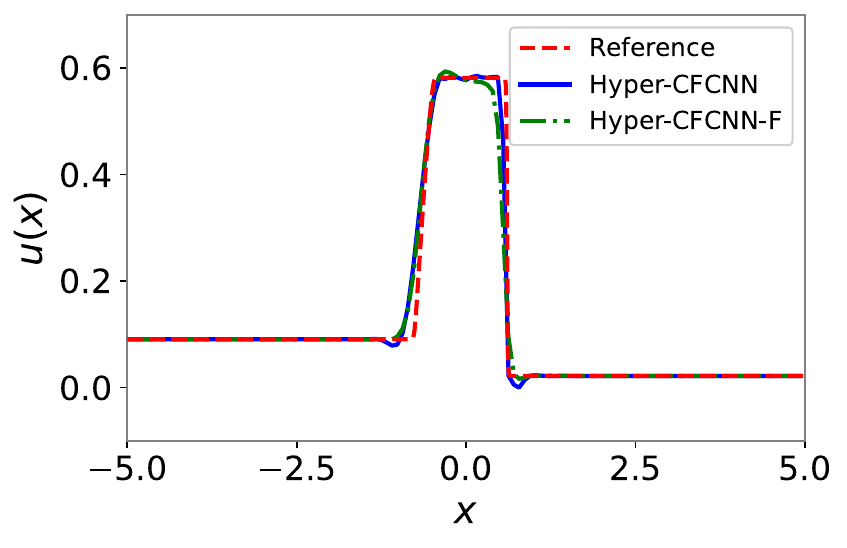}
\caption{$hv$, $T=1.0$, $N=128$}
\label{fig:sw_hv_T1_N128}
\end{subfigure}
\hfill
\begin{subfigure}{0.31\textwidth}
\includegraphics[width=\textwidth]{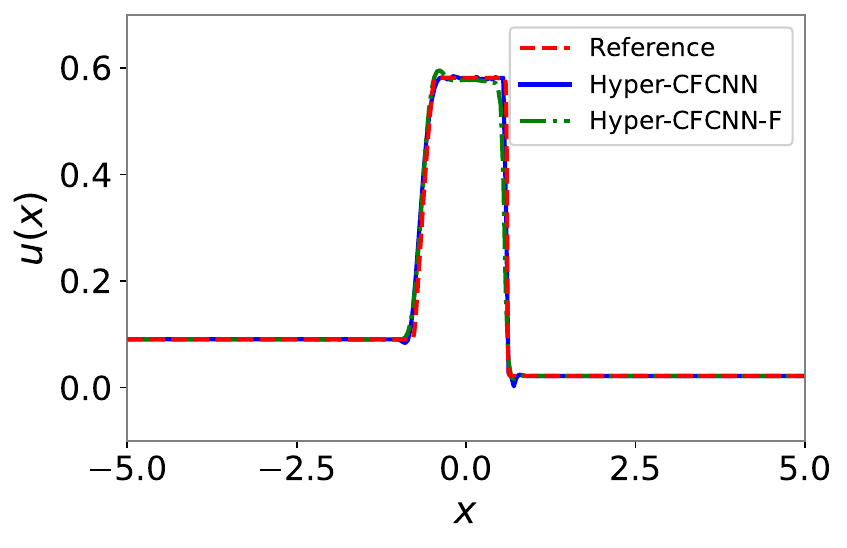}
\caption{$hv$, $T=1.0$, $N=256$}
\label{fig:sw_hv_T1_N256}
\end{subfigure}

\caption{H--CFCNN and H--CFCNN--F versus the WENO5
reference for the shallow-water example at \(T=0.5\) and \(T=1.0\) on \(N=64,128,256\).}
\label{fig:sw_convergence}
\end{figure}

\paragraph{Accuracy, conservation, and refinement}
Figure~\ref{fig:sw_convergence} compares H--CFCNN, H--CFCNN--F, and the WENO5
reference at \(T=0.5\) and \(T=1.0\) on \(N=64,128,256\). Both learned schemes
track the reference closely, and the discrepancy decreases as the mesh is
refined. Small oscillations visible on the coarse meshes are substantially
reduced on the finest grid.

\begin{figure}[htbp]
\centering
  \begin{subfigure}{0.48\textwidth}
    \includegraphics[width=\textwidth]{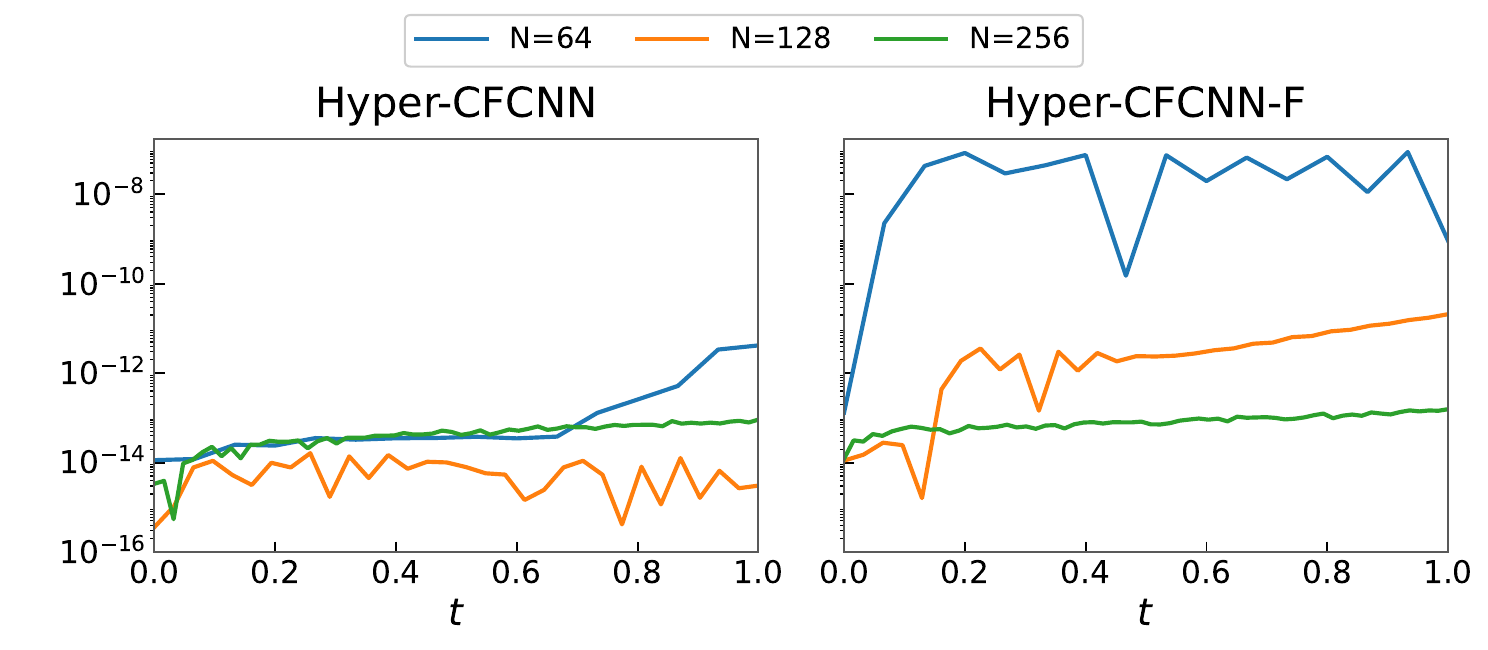}
    \caption{\(C(h)\)}
    \label{fig:shallow-mass_h}
  \end{subfigure}
  \hfill
   \begin{subfigure}{0.48\textwidth}
    \includegraphics[width=\textwidth]{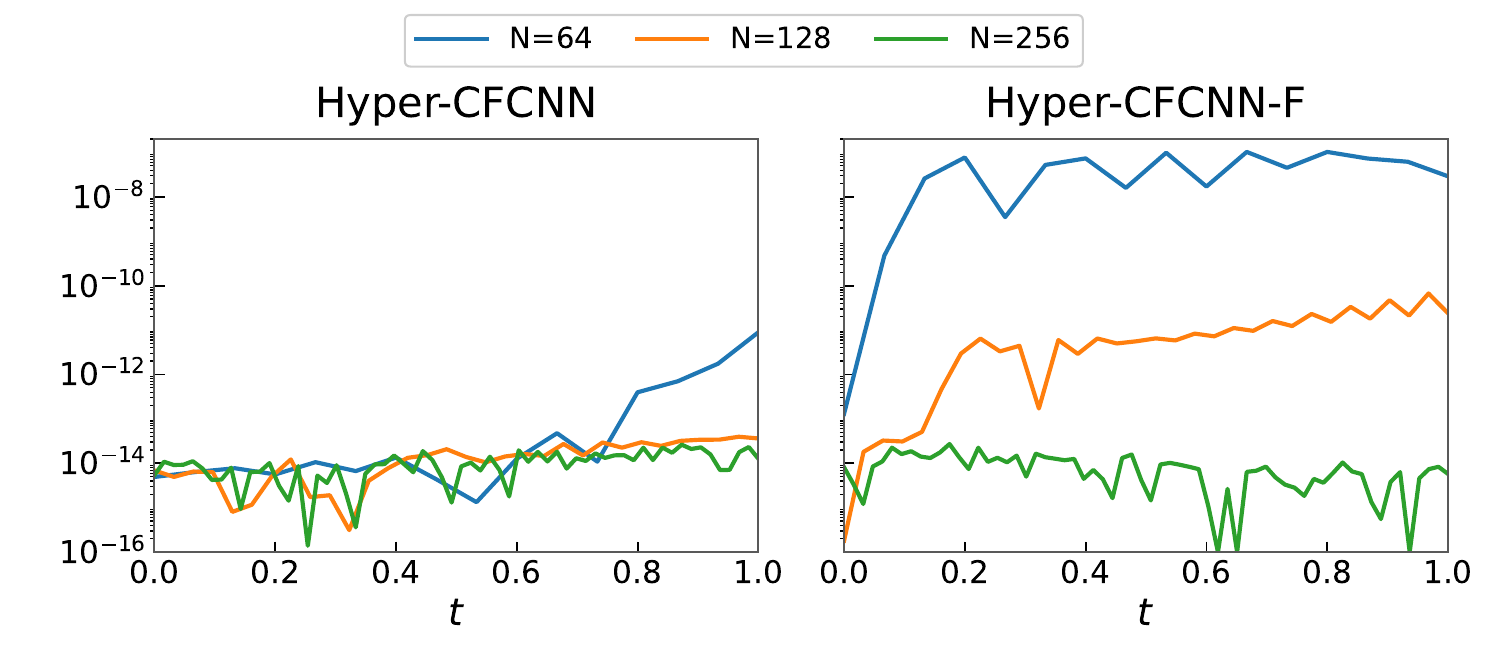}
    \caption{\(C(hv)\)}
    \label{fig:shallow-mass_hv}
  \end{subfigure}
  \hfill
  \caption{Conservation remainders \(C(h)\) and \(C(hv)\) for the shallow-water
example on \(N=64,128,256\).}
  \label{fig:sw_mass}
\end{figure}

Using \eqref{eq:disc_cons_remainder}, we track the remainders \(C(h)(t)\) and
\(C(hv)(t)\). Figure~\ref{fig:sw_mass} shows that both H--CFCNN and
H--CFCNN--F keep these remainders small, with values decreasing clearly under
mesh refinement. This behavior is consistent with the conservative
flux-difference form \eqref{eq:fv_net}.

Table~\ref{tab:sw_error} reports the MSE at \(T=1.0\) for
H--CFCNN together with the empirical refinement rates. The observed rates remain
below first order, which is expected in the presence of shocks and other
nonsmooth features, but the error decreases consistently under refinement.
Table~\ref{tab:sw_cost} reports the corresponding target-network parameter count
and wall-clock prediction time. As in the Burgers examples, the parameter count
and prediction cost increase linearly with the number of cells.

\begin{figure}[htbp]
  \centering
  \begin{subfigure}{0.23\textwidth}
    \includegraphics[width=\textwidth]{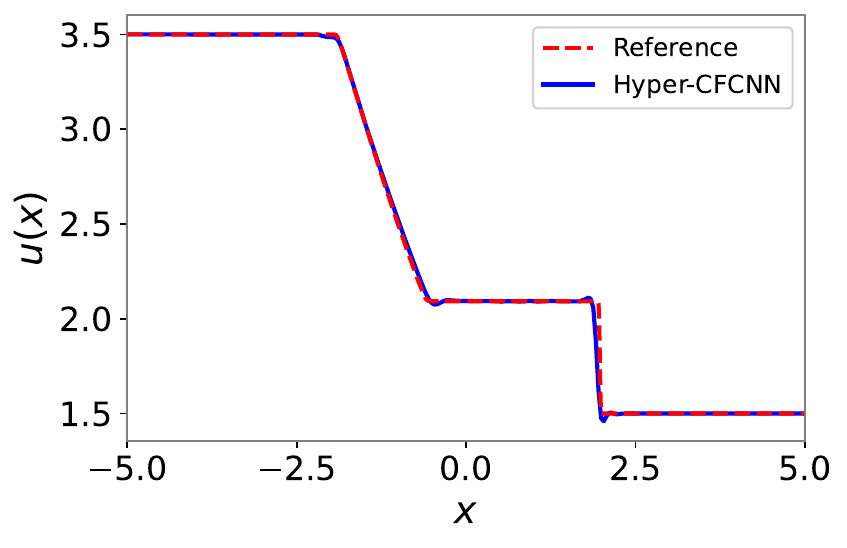}
    \caption{$h$ (init)}
    \label{fig:sw_init_h}
  \end{subfigure}
  \hfill
  \begin{subfigure}{0.23\textwidth}
    \includegraphics[width=\textwidth]{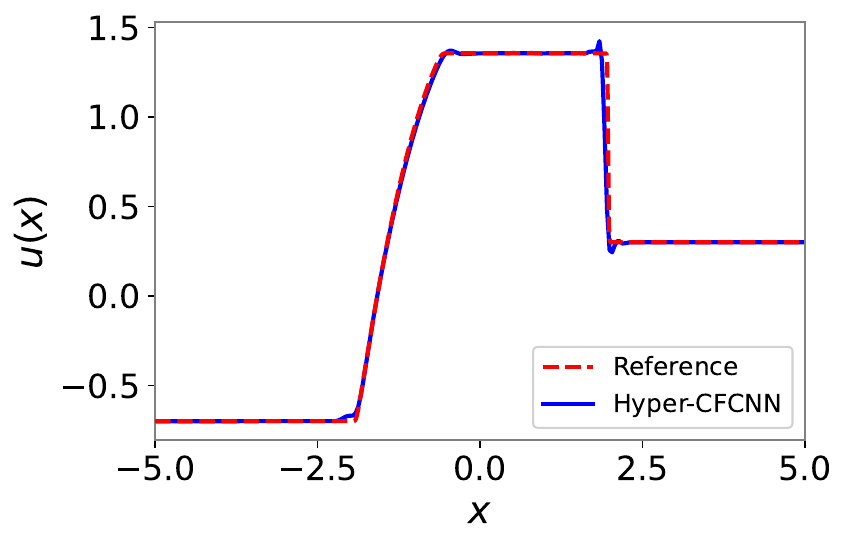}
    \caption{$hv$ (init)}
    \label{fig:sw_init_hv}
  \end{subfigure}
  \hfill
  \begin{subfigure}{0.23\textwidth}
    \includegraphics[width=\textwidth]{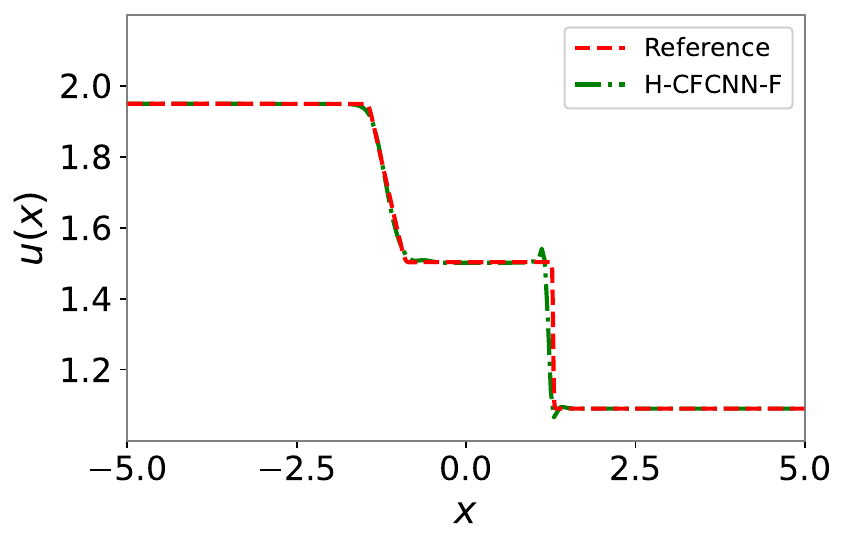}
    \caption{$h$ (mesh)}
    \label{fig:sw_mesh_h}
  \end{subfigure}
  \hfill
  \begin{subfigure}{0.23\textwidth}
    \includegraphics[width=\textwidth]{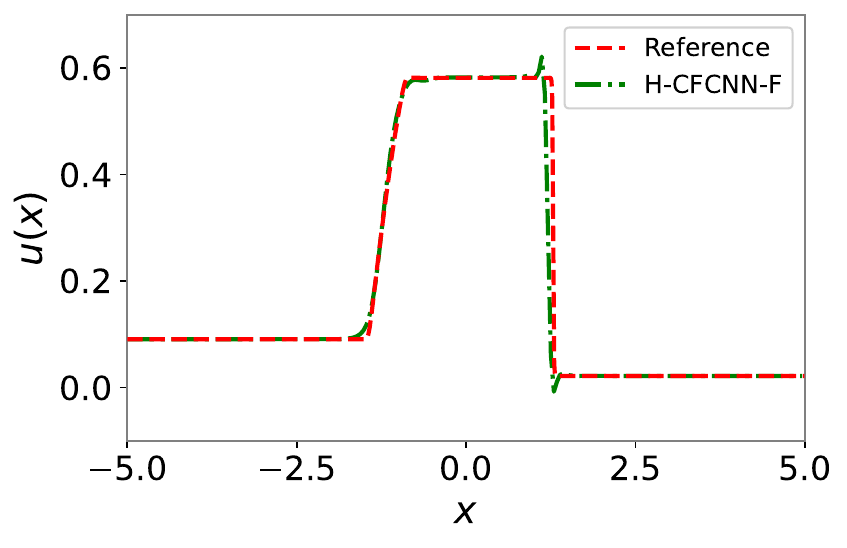}
    \caption{$hv$ (mesh)}
    \label{fig:sw_mesh_hv}
  \end{subfigure}
  \caption{H--CFCNN and H--CFCNN--F versus the reference for the shallow-water
example under initial-condition extrapolation on \(N=256\) at \(T=1.0\) and
mesh extrapolation on the unseen mesh \(N=224\) at \(T=1.0\).}
  \label{fig:sw_combined}
\end{figure}

\paragraph{Extrapolation in parameters and mesh resolution}
We consider two complementary extrapolation settings. For initial-condition
extrapolation, Figure~\ref{fig:sw_combined}(a) and
Figure~\ref{fig:sw_combined}(b) show a Riemann datum on \(N=256\) at \(T=1.0\)
with parameter values
\[
(h_\ell,h_r,v_\ell,v_r,x_0)=(3.5,\,1.5,\,-0.2,\,0.2,\,0.2),
\]
which lie outside the training ranges. In this case, H--CFCNN remains closer
to the reference for both \(h\) and \(hv\).

For mesh extrapolation, Figure~\ref{fig:sw_combined}(c) and
Figure~\ref{fig:sw_combined}(d) show results on the unseen mesh \(N=224\).
Here H--CFCNN--F produces accurate and stable predictions for both variables.
These results suggest that the known-flux and learned-flux variants have
different strengths in this example: the former is more robust under
extrapolation in the initial data, whereas the latter transfers well across
unseen mesh resolutions.

\begin{table}[!htbp]
\centering
\caption{MSE and empirical refinement rates for H--CFCNN on the shallow-water example at \(T=1.0\).}
\label{tab:sw_error}
\begin{tabular}{ccc}
\toprule
Mesh size \(N\) & MSE & Order \(p\) \\
\midrule
64  & \(5.473\times10^{-4}\) & --   \\
128 & \(1.830\times10^{-4}\) & 0.79 \\
256 & \(5.857\times10^{-5}\) & 0.82 \\
\bottomrule
\end{tabular}
\end{table}

\begin{table}[!htbp]
\centering
\caption{Target-network parameter count and wall-clock prediction time for H--CFCNN on the shallow-water example from \(t=0\) to \(t=1.0\).}
\label{tab:sw_cost}
\begin{tabular}{ccc}
\toprule
Mesh size \(N\) & Parameters & Time (s) \\
\midrule
64  & 2560  & 2.3108 \\
128 & 5120  & 4.6792 \\
256 & 10240 & 9.6471 \\
\bottomrule
\end{tabular}
\end{table}

\subsection{Euler equations}

We consider the one-dimensional compressible Euler equations for an ideal gas on
\(x\in[-5,5]\):
\begin{equation}\label{eq:Euler1D}
\begin{aligned}
\rho_t + (\rho u)_x &= 0,\\
(\rho u)_t + \bigl(\rho u^2 + p\bigr)_x &= 0,\\
E_t + \bigl(u(E+p)\bigr)_x &= 0,
\end{aligned}
\end{equation}
where \(\rho\) is the density, \(u\) the velocity, and \(E\) the total energy.
The pressure is given by
\begin{equation}\label{eq:EOS}
p=(\gamma-1)\Bigl(E-\tfrac{1}{2}\rho u^2\Bigr),\qquad \gamma=1.4.
\end{equation}

The initial condition is a randomized Shu--Osher profile, smoothed near the
right boundary:
\begin{equation}\label{eq:EulerIC}
(\rho,u,p)(x,0)=
\begin{cases}
(\rho_\ell,u_\ell,p_\ell), & x\le x_0,\\[3pt]
\bigl(1+\varepsilon\sin(5x),\,0,\,p_r\bigr), & x_0<x\le x_1,\\[3pt]
\bigl(1+\varepsilon\sin(5x)e^{-(x-x_1)^4},\,0,\,p_r\bigr), & x>x_1,
\end{cases}
\end{equation}
with nominal values
\(\hat\rho_\ell=3.857135\), \(\hat u_\ell=2.629369\),
\(\hat p_\ell=10.33333\), \(\hat p_r=1\), \(\hat x_0=-4\),
\(x_1=3.29867\), and \(\hat\varepsilon=0.2\). Training instances are obtained
by perturbing
\((\rho_\ell,u_\ell,p_\ell,p_r,\varepsilon,x_0)\)
uniformly within \(\pm 10\%\) of these nominal values.

\paragraph{Training and reference setup}
Reference trajectories are generated by a conservative finite-volume WENO5
solver with local Lax--Friedrichs flux and third-order SSP--RK time integration
on a uniform mesh with \(N=512\) and CFL\(=0.064\). Training is carried out on
the mesh family \(N\in\{64,128,192,256\}\), with time steps chosen to match the
same CFL across meshes, window length \(L=30\), and unroll depth \(K=4\). We
use \(N_{\mathrm{traj}}=200\) randomized initial conditions. In H--CFCNN--F,
the FluxNet parameters \(\boldsymbol{\theta}_{\mathrm{flux}}\) are trained
jointly with the hypernetwork; for this benchmark, \textsf{FluxNet} is taken to
be strictly local, with \(1\times1\) kernels. Unless stated otherwise, the accuracy, conservation, and refinement results
reported below use the fixed test parameters
\(
\rho_\ell=3.760352,
u_\ell=2.681251,
p_\ell=11.264806,
p_r=1.046399,
\varepsilon=0.186240,
x_0=-3.724815
\)
in \eqref{eq:EulerIC}.

\begin{figure}[htbp]
\centering

\begin{subfigure}{0.31\textwidth}
\includegraphics[width=\textwidth]{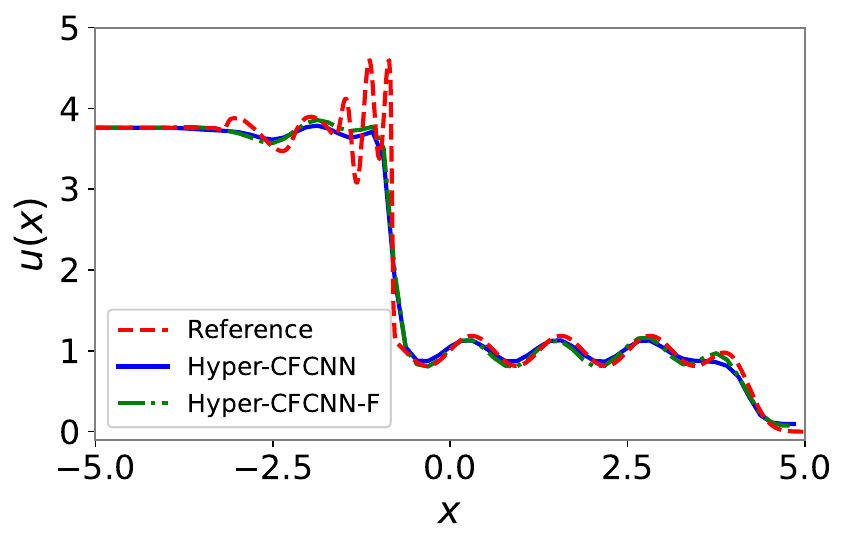}
\caption{$\rho$, $T=0.8$, $N=64$}
\label{fig:euler_rho_t08_N64}
\end{subfigure}
\hfill
\begin{subfigure}{0.31\textwidth}
\includegraphics[width=\textwidth]{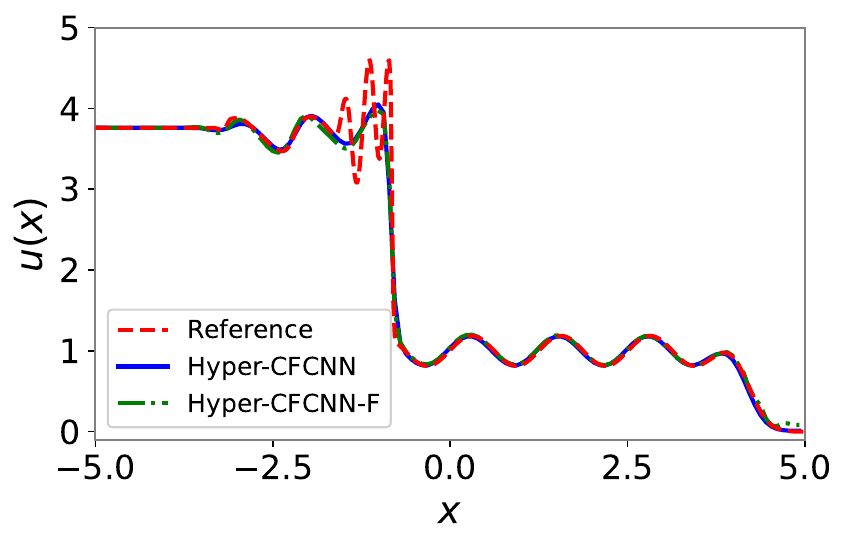}
\caption{$\rho$, $T=0.8$, $N=128$}
\label{fig:euler_rho_t08_N128}
\end{subfigure}
\hfill
\begin{subfigure}{0.31\textwidth}
\includegraphics[width=\textwidth]{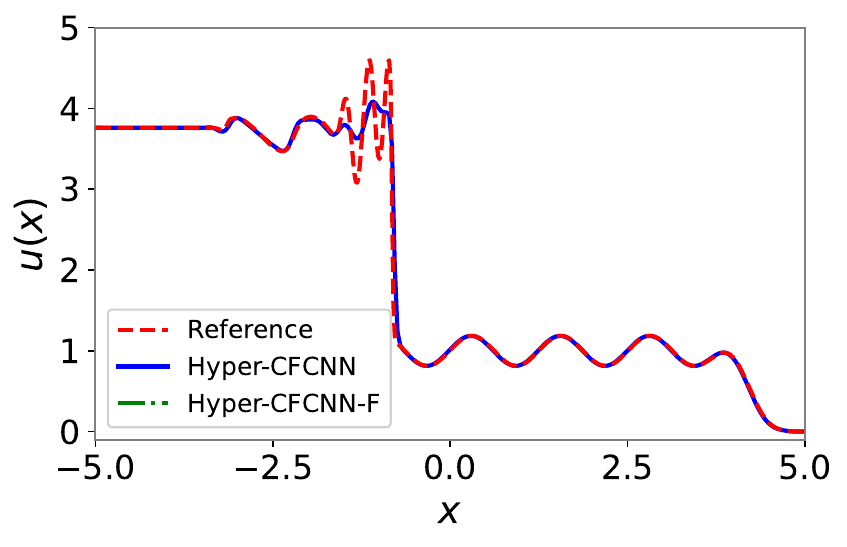}
\caption{$\rho$, $T=0.8$, $N=256$}
\label{fig:euler_rho_t08_N256}
\end{subfigure}

\vspace{0.8ex}

\begin{subfigure}{0.31\textwidth}
\includegraphics[width=\textwidth]{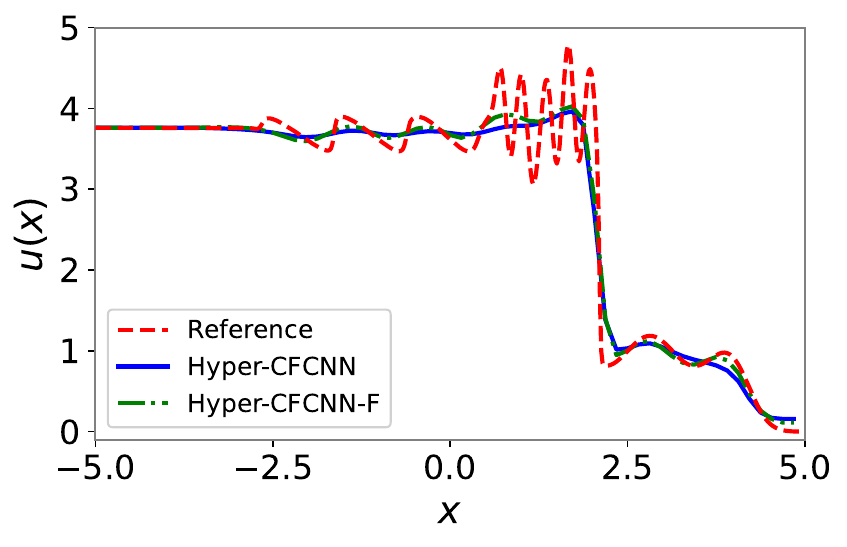}
\caption{$\rho$, $T=1.6$, $N=64$}
\label{fig:euler_rho_t16_N64}
\end{subfigure}
\hfill
\begin{subfigure}{0.31\textwidth}
\includegraphics[width=\textwidth]{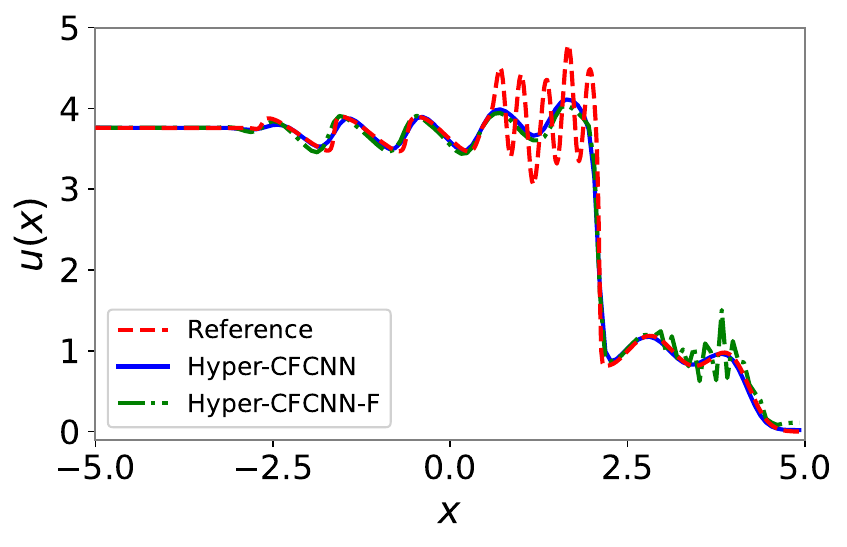}
\caption{$\rho$, $T=1.6$, $N=128$}
\label{fig:euler_rho_t16_N128}
\end{subfigure}
\hfill
\begin{subfigure}{0.31\textwidth}
\includegraphics[width=\textwidth]{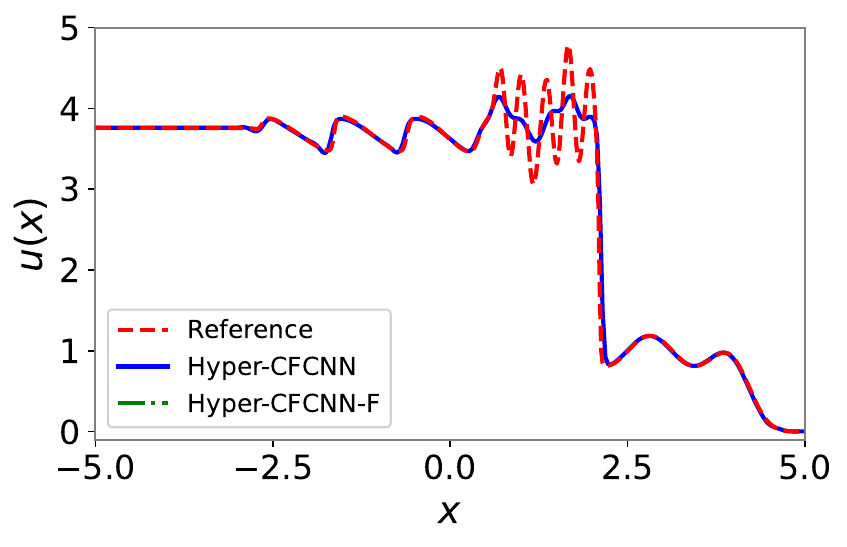}
\caption{$\rho$, $T=1.6$, $N=256$}
\label{fig:euler_rho_t16_N256}
\end{subfigure}

\vspace{0.8ex}

\begin{subfigure}{0.31\textwidth}
\includegraphics[width=\textwidth]{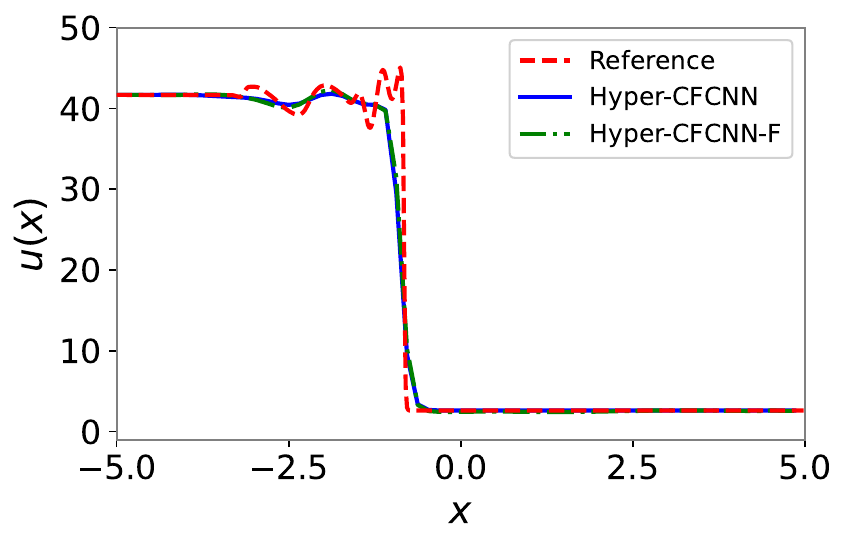}
\caption{$E$, $T=0.8$, $N=64$}
\label{fig:euler_E_t08_N64}
\end{subfigure}
\hfill
\begin{subfigure}{0.31\textwidth}
\includegraphics[width=\textwidth]{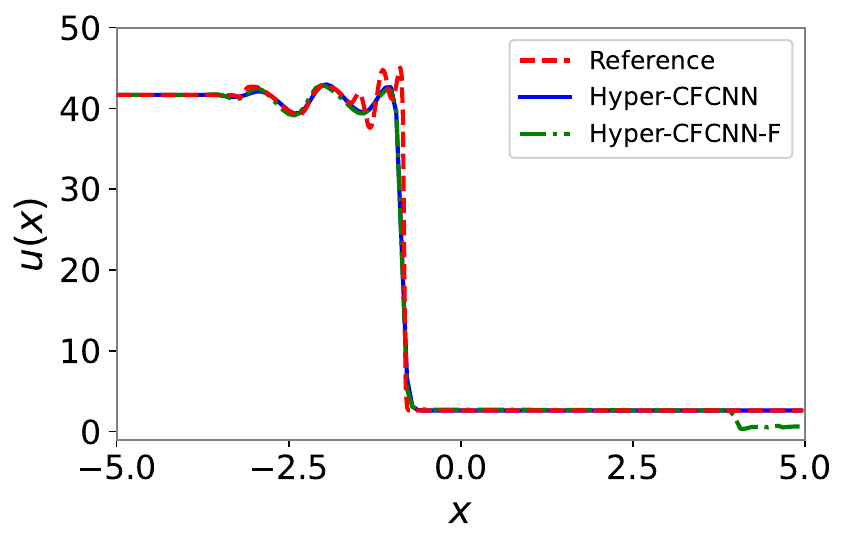}
\caption{$E$, $T=0.8$, $N=128$}
\label{fig:euler_E_t08_N128}
\end{subfigure}
\hfill
\begin{subfigure}{0.31\textwidth}
\includegraphics[width=\textwidth]{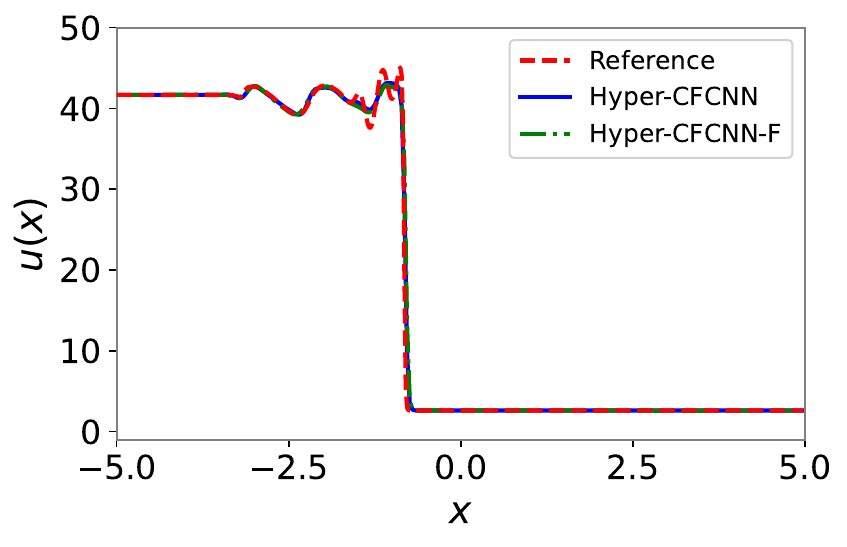}
\caption{$E$, $T=0.8$, $N=256$}
\label{fig:euler_E_t08_N256}
\end{subfigure}

\vspace{0.8ex}

\begin{subfigure}{0.31\textwidth}
\includegraphics[width=\textwidth]{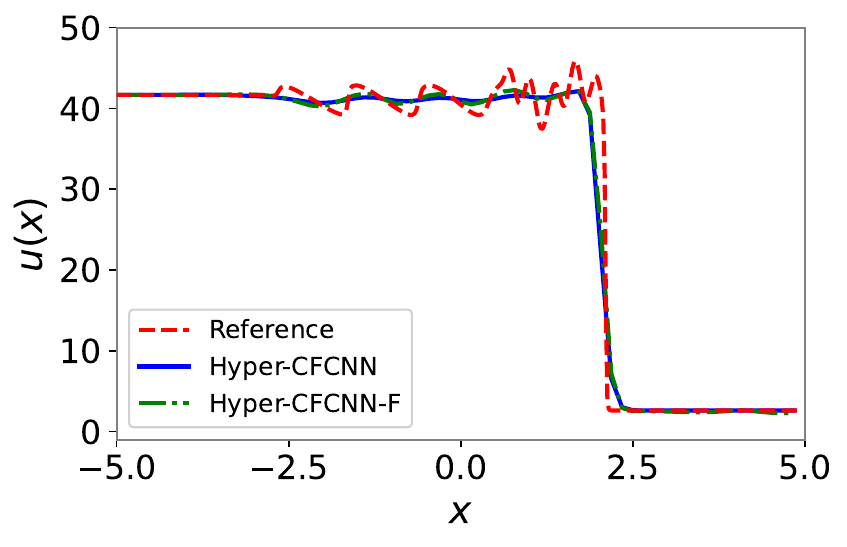}
\caption{$E$, $T=1.6$, $N=64$}
\label{fig:euler_E_t16_N64}
\end{subfigure}
\hfill
\begin{subfigure}{0.31\textwidth}
\includegraphics[width=\textwidth]{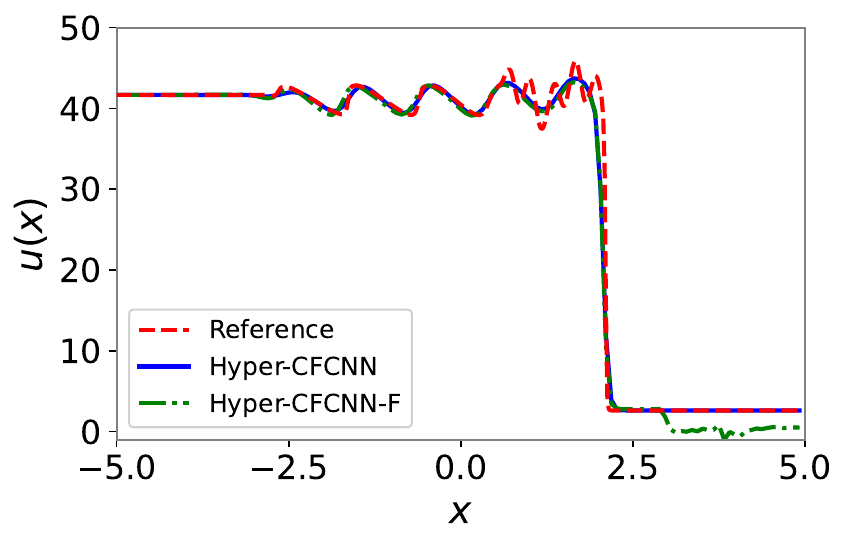}
\caption{$E$, $T=1.6$, $N=128$}
\label{fig:euler_E_t16_N128}
\end{subfigure}
\hfill
\begin{subfigure}{0.31\textwidth}
\includegraphics[width=\textwidth]{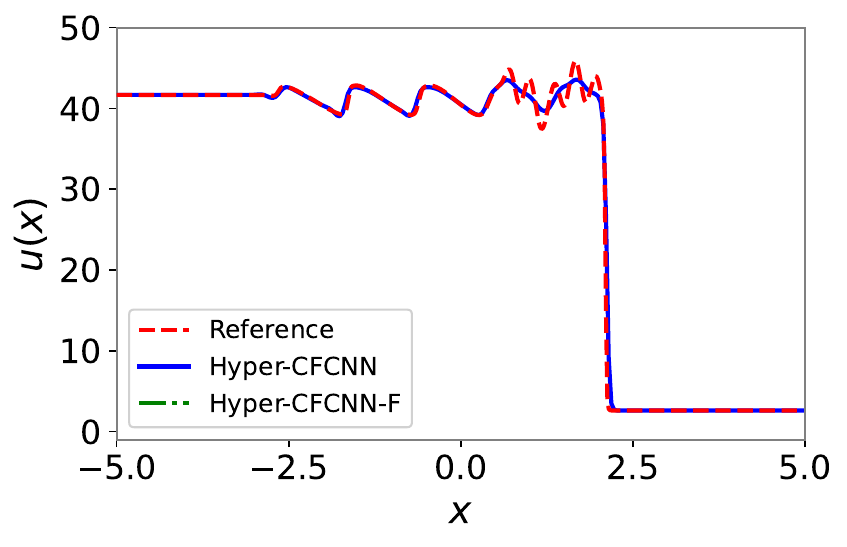}
\caption{$E$, $T=1.6$, $N=256$}
\label{fig:euler_E_t16_N256}
\end{subfigure}

\caption{H--CFCNN and H--CFCNN--F versus the WENO5 reference for the Euler
example, shown in the density \(\rho\) and total energy \(E\) at \(T=0.8\) and
\(T=1.6\) on \(N=64,128,256\).}
\label{fig:euler_rho_E}
\end{figure}

\paragraph{Accuracy, conservation, and refinement}
Figure~\ref{fig:euler_rho_E} compares the density \(\rho\) and total energy
\(E\) at \(T=0.8\) and \(T=1.6\) for \(N=64,128,256\). On the coarsest mesh,
both H--CFCNN and H--CFCNN--F capture the shock location and the downstream
oscillatory structure with visible phase and amplitude errors relative to the
reference. These discrepancies decrease under mesh refinement and are visually
small on \(N=256\). On the coarser meshes, H--CFCNN--F is slightly more
diffusive in the post-shock oscillatory region, but this difference becomes
less pronounced as the mesh is refined.

\begin{figure}[htbp]
  \centering
  \begin{subfigure}{0.48\textwidth}
    \centering
    \includegraphics[width=\textwidth]{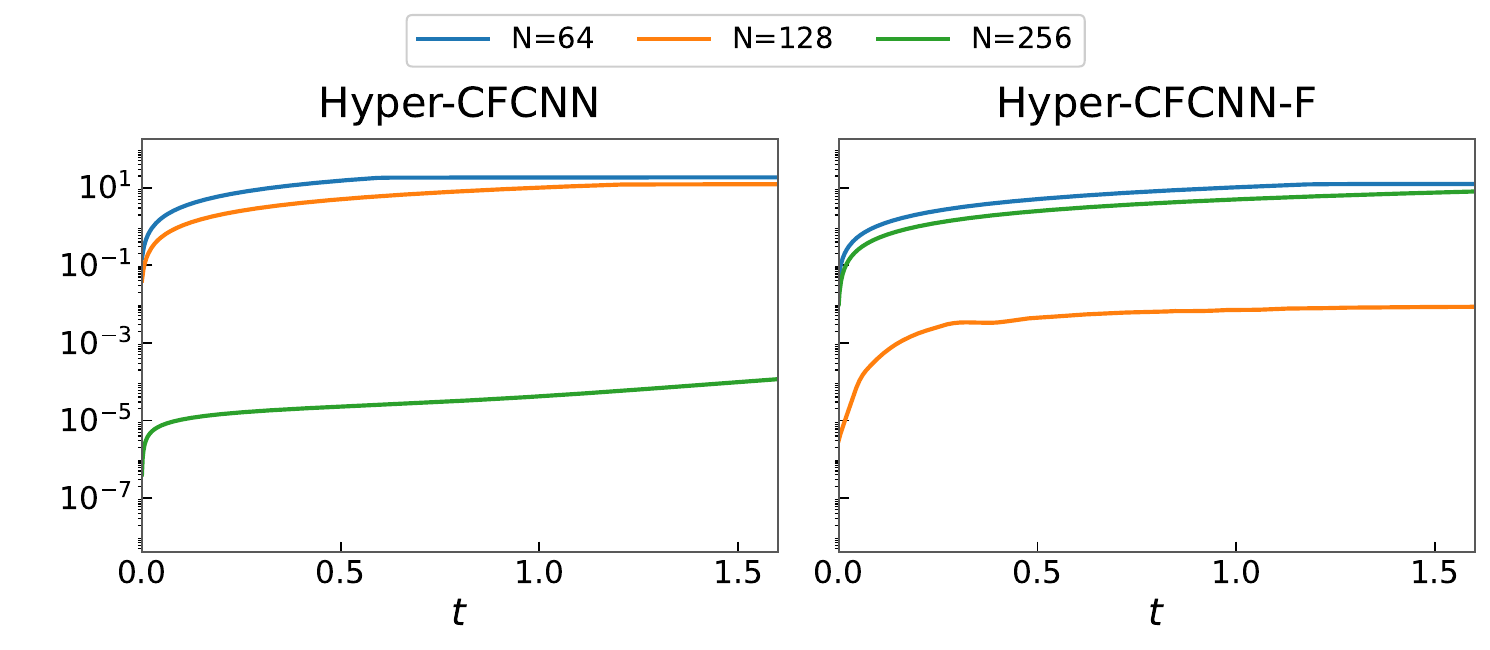}
    \caption{\(C(\rho)\)}
    \label{fig:euler_rho_cons}
  \end{subfigure}
  \hfill
  \begin{subfigure}{0.48\textwidth}
    \centering
    \includegraphics[width=\textwidth]{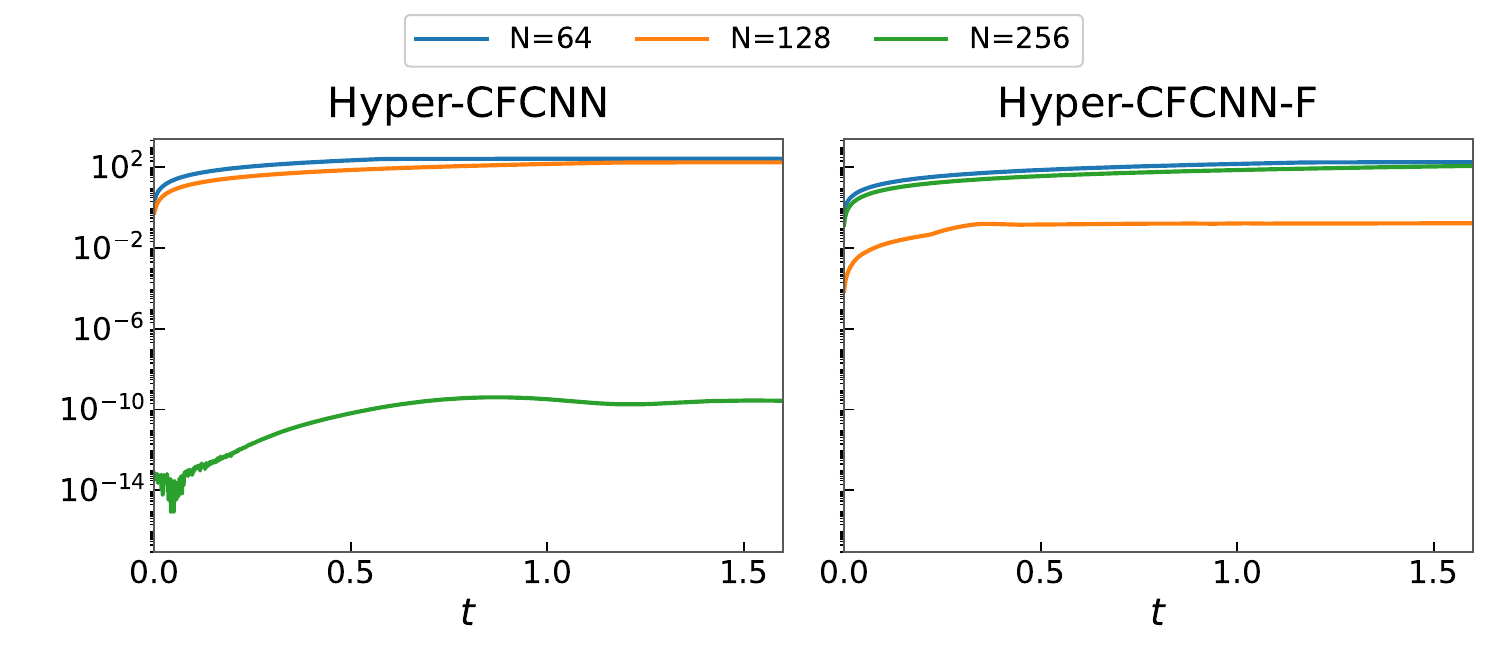}
    \caption{\(C(E)\) }
    \label{fig:euler_E_cons}
  \end{subfigure}

  \caption{Conservation remainders \(C(\rho)\) and \(C(E)\) for H--CFCNN and
H--CFCNN--F on \(N=64,128,256\) for the Euler example.}
  \label{fig:euler_conserve}
\end{figure}

Figure~\ref{fig:euler_conserve}(a) and Figure~\ref{fig:euler_conserve}(b)
report the conservation remainders \(C(\rho)(t)\) and \(C(E)(t)\) up to
\(T=1.6\). For H--CFCNN, the drift is substantially smaller on the finest mesh
than on the coarser meshes, and for total energy it is close to roundoff on
\(N=256\). For H--CFCNN--F, the remainders are larger and their dependence on
mesh resolution is less uniform, especially for total energy. This behavior is
consistent with the fact that the learned flux does not explicitly enforce the
algebraic structure satisfied by the analytical flux, even though the
conservative flux-difference form keeps the drift bounded over the rollout.

Table~\ref{tab:euler_cost} reports the target-network parameter count and the
wall-clock prediction time required to advance the solution from \(t=0\) to
\(t=1.6\) for H--CFCNN on representative meshes. As in the Burgers and
shallow-water examples, both quantities grow approximately linearly with the
number of cells.

\begin{figure}[!htbp]
  \centering
  \begin{subfigure}{0.48\textwidth}
    \centering
    \includegraphics[width=\textwidth]{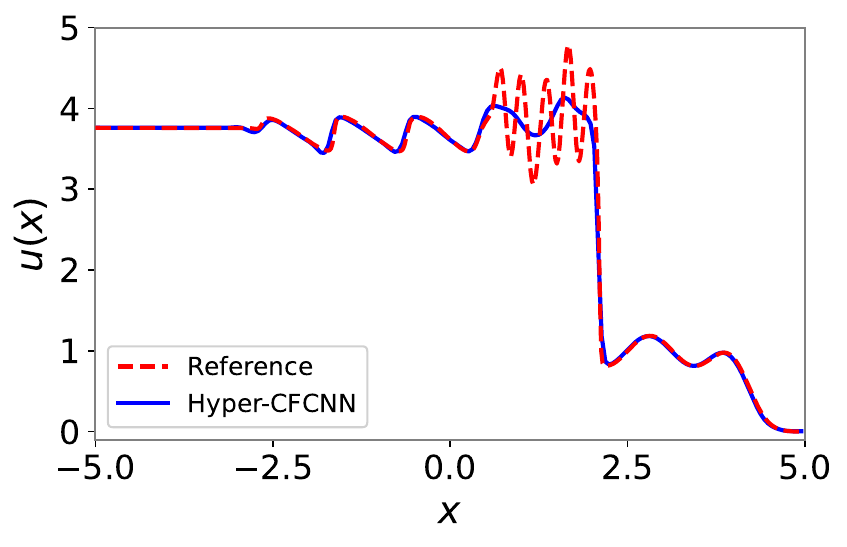}
    \caption{Density \(\rho\) on the unseen mesh \(N=224\) at \(T=1.6\).}
    \label{fig:euler_mesh_extra_rho}
  \end{subfigure}
  \hfill
  \begin{subfigure}{0.48\textwidth}
    \centering
    \includegraphics[width=\textwidth]{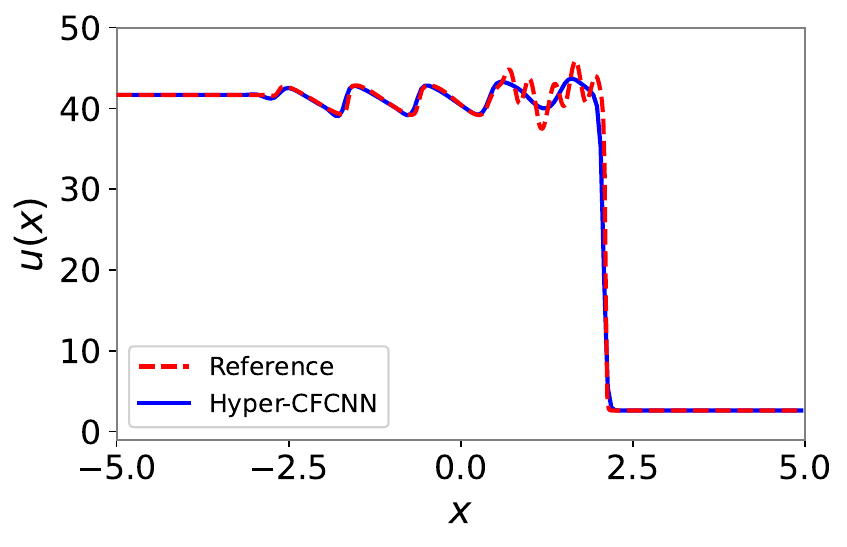}
    \caption{Total energy \(E\) on the unseen mesh \(N=224\) at \(T=1.6\).}
    \label{fig:euler_mesh_extra_E}
  \end{subfigure}
 \caption{H--CFCNN versus the reference for the Euler example on the unseen mesh
\(N=224\) at \(T=1.6\), shown in the density \(\rho\) and total energy \(E\).}
  \label{fig:euler_mesh_extra}
\end{figure}

\paragraph{Mesh transfer}
Figure~\ref{fig:euler_mesh_extra} examines transfer to the unseen mesh
\(N=224\), which lies inside the trained resolution interval but is not used
during training. H--CFCNN reproduces the main shock location, the post-shock
oscillatory structure, and the downstream state with good fidelity. In the
high-frequency region behind the shock, the prediction is slightly smoother than
the reference, especially in total energy, but the overall agreement remains
good.

The Euler example is more demanding than the scalar tests because it involves a
system, strong gradients, and a longer oscillatory post-shock region. In this
setting, H--CFCNN is more accurate and more conservative than H--CFCNN--F,
especially on coarse meshes and in total energy. The learned-flux variant
nevertheless captures the main shock dynamics and improves under mesh
refinement, which supports its use when the analytical flux is unavailable.

\begin{table}[!htbp]
\centering
\caption{Target-network parameter count and wall-clock prediction time for H--CFCNN on the Euler example from \(t=0\) to \(t=1.6\).}
\label{tab:euler_cost}
\begin{tabular}{ccc}
\toprule
Mesh size \(N\) & Parameters & Time (s) \\
\midrule
64  & 2560  & 39.22  \\
128 & 5120  & 78.79  \\
256 & 10240 & 156.69 \\
\bottomrule
\end{tabular}
\end{table}

\section{Conclusion}\label{sec:con}

We have studied a conservative data-driven finite-volume framework for
one-dimensional nonlinear conservation laws in which only selected components of
the discretization are learned. In the known-flux setting, Hyper--CFCNN learns
the nonlinear WENO5 weight map while retaining the classical reconstruction,
conservative flux-difference form, and SSP--RK time integration. In the
unknown-flux setting, Hyper--CFCNN--F further replaces the analytical flux by a
learned interface-flux model within the same finite-volume update structure. The numerical results on Burgers, shallow-water, and one-dimensional Euler
examples show that both variants improve under mesh refinement and can be
applied across multiple mesh resolutions without retraining. Overall,
H--CFCNN provides the stronger accuracy and conservation behavior, especially on
the Euler example, where the analytical flux remains important for controlling
drift in conserved quantities. H--CFCNN--F is generally less accurate and less
conservative on coarse meshes, but it still captures the main solution
structure and improves consistently under refinement. The experiments also
suggest that the two variants have different strengths: the known-flux
formulation is more robust under extrapolation in problem parameters, whereas
the learned-flux formulation transfers well across mesh resolutions.

Several directions remain open. On the theoretical side, it would be useful to
identify conditions under which the learned reconstruction recovers the
smooth-region accuracy of WENO5, and to analyze stability-related properties
such as entropy stability and positivity preservation, especially for systems.
On the numerical side, broader comparisons with stronger analytical baselines,
including WENO-Z and mapped WENO, together with systematic cost--accuracy
studies, would help clarify the regimes in which the learned variants are most
effective. Extending the framework to multi-dimensional problems and to a wider
class of boundary conditions is also of interest.

\section*{Acknowledgments}

\bibliographystyle{siamplain}
\bibliography{main}
\end{document}